\documentclass[10pt]{article}
\usepackage[T1]{fontenc}
\usepackage[utf8]{inputenc} % pour gÃÂ'Â''ÂÃÂ'Â''ÂÃÂ'Â''ÂÃÂ'Â''ÂÃÂ'Â''ÂÃÂ'Â''ÂÃÂ'Â''ÂÃÂ'Â''ÂÃÂ'Â''ÂÃÂ'Â''ÂÃÂ'Â''ÂÃÂ'Â''ÂÃÂ'Â''ÂÃÂ'Â''ÂÃÂ'Â''ÂÃÂ'Â''Â¡?rer les caractÃÂ'Â''ÂÃÂ'Â''ÂÃÂ'Â''ÂÃÂ'Â''ÂÃÂ'Â''ÂÃÂ'Â''ÂÃÂ'Â''ÂÃÂ'Â''ÂÃÂ'Â''ÂÃÂ'Â''ÂÃÂ'Â''ÂÃÂ'Â''ÂÃÂ'Â''ÂÃÂ'Â''ÂÃÂ'Â''ÂÃÂ'Â''Â¡?res accentuÃÂ'Â''ÂÃÂ'Â''ÂÃÂ'Â''ÂÃÂ'Â''ÂÃÂ'Â''ÂÃÂ'Â''ÂÃÂ'Â''ÂÃÂ'Â''ÂÃÂ'Â''ÂÃÂ'Â''ÂÃÂ'Â''ÂÃÂ'Â''ÂÃÂ'Â''ÂÃÂ'Â''ÂÃÂ'Â''ÂÃÂ'Â''Â¡?es

\usepackage[english]{babel}
\usepackage{amsmath}
\usepackage{amssymb,mathrsfs}
\usepackage{enumitem} 
\usepackage{amsthm}
\usepackage{titlesec}
\usepackage[all]{xy}
\usepackage[filecolor=green]{hyperref}

\usepackage[a4paper]{geometry}
\def\l@paragraph{\@tocline{3}{0pt}{1pc}{9pc}{}}

\newcommand\M[1]{\mathscr{#1}}

\newcommand\N{\mathbb{N}}

\newcommand\K[1]{\mathbb{K}{#1}}
\newcommand\red[1]{\text{Red}\left(#1\right)}

\newcommand\nred[1]{\text{Nred}\left(#1\right)}
\newcommand\im[1]{\text{im}\left(#1\right)}
\newcommand\redd[1]{\emph{Red}\left(#1\right)}
\newcommand\nredd[1]{\emph{Nred}\left(#1\right)}

\newcommand\obsred[1]{\text{Obs}^{#1}}
\newcommand\obsredd[1]{\emph{Obs}^{#1}}

\newcommand\F[1]{\underset{#1}{\longrightarrow}}

\newcommand\X[1]{X^{\left({#1}\right)}}

\newcommand\EV[1]{\left\langle #1\right\rangle}

\newcommand\id[1]{\text{Id}_{#1}}
\newcommand\idd[1]{\emph{Id}_{#1}}

\newcommand\lm[1]{\text{lm}\left(#1\right)}
\newcommand\lmm[1]{\emph{lm}\left(#1\right)}
\newcommand\lgen[1]{\text{lg}\left(#1\right)}

\newcommand\lc[1]{\text{lc}\left(#1\right)}
\newcommand{\G}{Gr\"obner}
\newcommand{\RO}{\textbf{RO}\left(G,<\right)}

\titleformat{\subsubsection}[runin]
{\normalfont\bfseries}
{\thesubsubsection.}{.5em}{}[.]
\titlespacing{\subparagraph}
{\parindent}{1.5ex plus .1ex minus .2ex}{1.5pt}
\begin{document}

\title{Reduction Operators and Completion of Rewriting Systems}
\author{CYRILLE CHENAVIER}
\date{}
\maketitle

\begin{abstract}

We propose a functional description of rewriting systems where reduction rules are represented by linear maps called reduction operators. We show that reduction operators admit a lattice structure. Using this structure we define the notions of confluence and of Church-Rosser property. We show that these notions are equivalent. We give an algebraic formulation of completion and show that such a completion exists using the lattice structure. We interpret the confluence for reduction operators in terms of \G\ bases. Finally, we introduce generalised reduction operators relative to non totally ordered sets.

\end{abstract}

\tableofcontents

\section{Introduction}

Convergent rewriting systems are confluent and terminating rewriting systems. They appear in rewriting theory to solve decision problems such as the word problem or the ideal membership problem. Completion algorithms were introduced to compute convergent rewriting systems: the Knuth-Bendix completion algorithm~\cite{MR0255472} for term rewriting~\cite{MR1629216} and string rewriting~\cite{MR1215932} or the Buchberger algorithm for \G\ bases \cite{buchberger1965algorithmus, MR2547481, MR0506423, MR506890} of commutative algebras \cite{buchberger1965algorithmus, MR893184} or associative algebras \cite{MR1299371}. In this paper, we propose an algebraic approach to completion: we formulate it algebraically and show that it can be obtained with an algebraic construction. 

We use the functional point of view considered by Berger~\cite{MR1608711} for rewriting on non-commutative polynomials. The latter are linear combinations of words. In this introduction, we first explain how the functional approach to string rewriting systems works. In the second part, we introduce \emph{reduction operators} and formulate the confluence and the completion with those. We also make explicit the link between reduction operators and rewriting on non-commutative polynomials, which gives us our algebraic constructions.\newpage

\begin{center}

\textbf{\large A Functional Approach to String Rewriting and \G\ bases}

\end{center}

\paragraph{Confluence for String Rewriting Systems.}

For string rewriting systems, the method consists in considering an idempotent application modelling the rewrite rules. This method works for \emph{semi-reduced} string rewriting systems, that is the systems such that
\begin{enumerate}
\item the left-hand sides of its rewrite rules are pairwise distinct,
\item no right-hand side of its rules is the left-hand side of another one.
\end{enumerate}
For instance, the string rewriting system with alphabet $\{x,y\}$ and with one rewrite rule $yy\F{}yx$ is semi-reduced. 

Given a string rewriting system $\EV{X\mid R}$ with alphabet $X$ and set of rewrite rules $R$, we denote by $X^*$ the set of words over $X$. Our algebraic constructions require that $\EV{X\mid R}$ is equipped with a total termination order $<$, that is, a terminating order on words such that every left-hand side of a rewrite rule is greater than the corresponding right-hand side. In Theorem~\ref{Unicity of normalised echelon basis}, we show, using this order, that $\EV{X\mid R}$ can be transformed into a unique semi-reduced string rewriting system, so that we may assume that it has this property. The application modelling its rules is the map $S:X^*\F{}X^*$ defined by
\begin{enumerate}
\item $S(l(\alpha))=r(\alpha)$ for every $\alpha\in R$ with left-hand side $l(\alpha)$ and right-hand side $r(\alpha)$,
\item $S(w)=w$ if no element of $R$ has left-hand side $w$.
\end{enumerate}
The application associated to our example maps $yy$ to $yx$ and fixes all other words.

The order $<$ guarantees that $\EV{X\mid R}$ terminates. Thus, it is sufficient to study whether it is confluent or not to know if it is convergent. In order to obtain the functional formulation of confluence, we consider the \emph{extensions of S}, that is, the applications $S_{p,q}$ defined for every pair of integers $\left(p,q\right)$ by
\begin{enumerate}
\item $S_{p,q}(w)=w_1r(\alpha)w_2$, if there exist words $w_1,\ w_2$ of length $p$ and $q$, respectively, and $\alpha\in R$, such that $w$ is equal to $w_1l(\alpha)w_2$,
\item $S_{p,q}(w)=w$, otherwise.
\end{enumerate}
In the previous example, $S_{0,1}$ maps $yyx$ to $yxx$ and $yyy$ to $yxy$, and $S_{1,0}$ maps $xyy$ to $xyx$ and $yyy$ to $yyx$. These applications enable us to characterise the normal forms for $\EV{X\mid R}$: a normal form is a word whose every sub-word is fixed by $S$, that is, the normal forms are the words fixed by all the extensions of $S$.

Given a word $w$, we denote by $[w]$ the class of $w$ for the equivalence relation induced by $R$. The order $<$ being total and well-founded, $[w]$ admits a smallest element. Let $M$ be the application from $X^*$ to itself mapping a word to this minimum. A word $w$ fixed by all the extensions of $S$ but which is not fixed by $M$ is called an \emph{obstruction} of $\EV{X\mid R}$. In other words, an obstruction is a normal form which is not minimal in its equivalence class. Hence, the set of obstructions is empty if and only if each equivalence class contains exactly one normal form. Moreover, recall that a terminating rewriting system is confluent if and only if every element admit exactly one normal form (see for instance~\cite[Section 2.1]{MR1629216}). Thus, we obtain the following functional characterisation of confluence: $\EV{X\mid R}$ is confluent if and only if the set of obstructions is empty. Considering our example, we deduce from the diagram

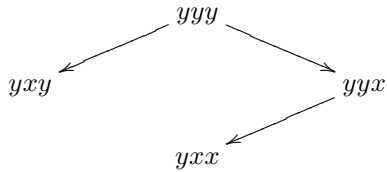
\[
\xymatrix @C = 4em @R = 1.5em{
&
yyy
\ar@1 [rd] 
\ar@1 [ld] 
& \\
yxy
&
&
yyx
\ar@1 [ld] 
\\
&
yxx
&
}
\]
that $yxx$ and $yxy$ belong to the same equivalence class. Moreover, each sub-word of $yxx$ and $yxy$ is fixed by $S$. Given a total order on the set of words, $yxx$ is either strictly smaller than $yxy$ or is strictly greater than $yxy$. In the first case, $yxy$ is an obstruction while in the second case $yxx$ is.

\paragraph{\G\ Bases and Homogeneous Rewriting Systems.}

Given a set $X$, we denote by $\K{X^*}$ the vector space spanned by $X^*$ over a commutative field $\K{}$: the non-zero elements of this vector space are the finite linear combinations of words with coefficients in $\K{}$. Let $<$ be a well-founded total order on $X^*$. Consider a set of rewrite rules $R$ on $\K{X^*}$ oriented with respect to $<$: for every $\alpha\in R$, $l(\alpha)$ is a word and is strictly greater than every word occurring in the decomposition of $r(\alpha)$ with respect to the basis $X^*$. We say that $R$ is a \emph{\G\ basis} when it induces a convergent rewriting system. The set $X^*$ is naturally embedded into $\K{X^*}$, so that a string rewriting system $\EV{X\mid R}$ induces a unique rewriting system on $\K{X^*}$. Moreover, $\EV{X\mid R}$ is convergent if and only if $R$, seen as a set of rewrite rules on $\K{X^*}$, is a \G\ basis.

The functional characterisation of the confluence for string rewriting systems extends into a functional characterisation of \G\ bases. The notion of semi-reduced string rewriting system is extended as follows: a rewriting system on $\K{X^*}$ with set of rewrite rules $R$ is said to be \emph{semi-reduced} if
\begin{enumerate}
\item the left-hand sides of the elements of $R$ are pairwise distinct words,
\item for every $\alpha,\ \beta\in R$, $l(\alpha)$ does not occur to the decomposition of $r(\beta)$ with respect to the basis $X^*$.
\end{enumerate}
As for string rewriting systems, Theorem~\ref{Unicity of normalised echelon basis} enables us to conclude that every rewriting system on $\K{X^*}$ can be transformed into a unique semi-reduced one. From now on, we assume that the rewriting system induced by $X$ and $R$ is semi-reduced.

The application mapping every left-hand side of a rewrite rule to its right-hand side induces an idempotent linear endomorphism $\overline{S}$ of $\K{X^*}$. For every pair of integers $(p,q)$, we consider the extension of $\overline{S}$
\[\overline{S}_{p,q}=\id{\left(\K{X}\right)^{\otimes p}}\otimes\overline{S}\otimes\id{\left(\K{X}\right)^{\otimes q}},\]
where for every integer $n$, $\left(\K{X}\right)^{\otimes n}$ denotes the $n$-fold tensor product of $\K{X}$. The operator $\overline{S}_{p,q}$ is the linear version of the function $S_{p,q}$ defined in the previous section: it is defined on the basis $X^*$ of $\K{X^*}$ in the following way
\begin{enumerate}
\item $\overline{S}_{p,q}(w)=w_1r(\alpha)w_2$, if there exist words $w_1,\ w_2$ of length $p$ and $q$, respectively, and $\alpha\in R$, such that $w$ is equal to $w_1l(\alpha)w_2$,
\item $\overline{S}_{p,q}(w)=w$, otherwise.
\end{enumerate}

Let $\overline{M}$ be the endomorphism of $\K{X^*}$ mapping every element $f\in\K{X^*}$ to the smallest element of $[f]$ for the natural multi-set order on $\K{X^*}$ induced by $<$. When one has a string rewriting system, $\overline{S}$, $\overline{S}_{p,q}$ and $\overline{M}$ are the linear endomorphisms of $\K{X^*}$ extending $S$, $S_{p,q}$ and $M$ defined in the previous paragraph, respectively. The reasoning we made also works in this context, so that we obtain: $R$ is a \G\ basis if and only if the set of obstructions is empty.

In ~\cite{MR1608711}, Berger considered this functional characterisation of \G\ bases to study \emph{finitely homogeneous rewriting systems}, that is, the systems such that $X$ is finite and there exists an integer $N$ such that the left-hand side and the right-hand side of every element of $R$ are linear combinations of words of length $N$. In particular, the endomorphism $\overline{S}$ associated to such a system induces an endomorphism of the vector space $\K{X}^{\otimes N}$ spanned by the set $\X{N}$ of words of length $N$. More generally, for every integer $n$, the extensions of $\overline{S}$ induce endomorphisms of the vector spaces spanned by the finite sets $\X{n}$. We denote by $F_n$ the set of endomorphisms of $\K{X}^{\otimes n}$ obtained with the restrictions of the extensions of $\overline{S}$. Moreover, the rewrite rules being homogeneous, for every word $w$ of length $n$, $\overline{M}(w)$ is a linear combination of words of length $n$, so that $\overline{M}$ also induces endomorphisms $\overline{M}_n$ of $\K{X}^{\otimes n}$. Hence, the set of obstructions admits a filtration on the length: an obstruction of length $n$ is a word fixed by every element of $F_n$ but not fixed by $\overline{M}_n$. Using the fact that each set $\X{n}$ is finite, Berger proved that $\overline{M}_n$ can be obtained from $F_n$ by an algebraic construction and deduced from this an algebraic formulation of obstructions, and thus, of \G\ bases for homogeneous rewriting systems.

This formulation enables us to obtain various proofs of \emph{Koszulness}~\cite{MR1608711, MR1683270, MR1832913, MR3299599, MR3503238}. Koszulness has applications in various topics: representation theory, numbers theory, algebraic and non-commutative geometry, for instance. We refer the reader to~\cite{MR0265437, MR1832913} for the definition of Koszul algebras and to~\cite{MR2177131} for an inventory of references about their applications.

\paragraph{Confluence for Non-Homogeneous Rewriting Systems.}

Consider a set of rewrite rules on $\K{X^*}$. When this set is non-homogeneous, $\overline{M}$ does not induce endomorphisms of $\K{X}^{\otimes n}$, so that we cannot construct it by restrictions on finite-dimensional vector spaces. Our first contribution is to show that it can be constructed globally on $\K{X^*}$. This construction uses the notion of reduction operator which are generalisations of the endomorphisms associated to a rewriting system on $\K{X^*}$.

\begin{center}

\textbf{\large Our Results}

\end{center}

\paragraph{Reduction Operators: Lattice Structure and Confluence.}

Let $G$ be a set and let $<$ be a well-founded total order on $G$. Typically, when we consider homogeneous rewriting systems, $G$ designates the sets $\X{n}$ and when we consider non-homogeneous rewriting systems, $G$ is the set $X^*$. A reduction operator relative to $\left(G,<\right)$ is a linear endomorphism $T$ of $\K{G}$ such that
\begin{enumerate}
\item $T$ is idempotent,
\item for every $g\in G$, $T(g)$ is either equal to $g$ or is a linear combination of elements of $G$ strictly smaller than $g$ for $<$.
\end{enumerate}

The set of reduction operators, written $\RO$, admits a lattice structure. Indeed, the first result of the paper about reduction operators is Proposition~\ref{Bijection} which states that the map $T\F{}\ker(T)$ from $\RO$ to the set of subspaces of $\K{G}$ is a bijection. This result extends, with a different method, the one of Berger who obtained it when $G$ is finite. The set of subspaces of $\K{G}$ admits a lattice structure: the order is the inclusion, the lower bound is the intersection and the upper bound is the sum. Using this structure as well as the bijection induced by the kernel map, we deduce a lattice structure on $\RO$.

Given a subset $F$ of $\RO$, let $\wedge F$ be its lower bound. We denote by $\red{T}$ the set of elements of $G$ fixed by a reduction operator $T$. In Lemma~\ref{Inclusion of images}, we show that $\red{\wedge F}$ is included in the intersection of all the $\red{T}$ where $T$ belongs to $F$:
\begin{equation} \label{Definition of obstructions I}
\red{\wedge F}\ \subseteq\ \bigcap_{T\in F}\red{T}.
\end{equation}
The complement of the inclusion (\ref{Definition of obstructions I}) is written $\obsred{F}$. The set $F$ is said to be \emph{confluent} if $\obsred{F}$ is empty.

Let $X$ be a set and let $R$ be a set of rewrite rules on $\K{X^*}$, oriented with respect to a well-founded total order on $X^*$. The endomorphism $\overline{S}$ associated to $R$, and more generally all the extensions of $\overline{S}$, are reduction operators relative to $\left(X^*,<\right)$. Let $F$ be the set of the extensions of $\overline{S}$. In Section~\ref{Rewriting properties of reduction operators}, we show that $\wedge F$ maps every element $f$ of $\K{X^*}$ to the smallest element of $[f]$, that is, $\wedge F$ is the operator $\overline{M}$ and $\obsred{F}$ is the set of obstructions of $\EV{X\mid R}$. We obtain our characterisation of \G\ bases in terms of reduction operators: $R$ is a \G\ basis if and only if $F$ is confluent.

\paragraph{Completion.}

Given a rewriting system on a set of terms, words or non-commutative polynomials with a set of rewrite rules $R$, the Knuth-Bendix completion algorithm or the Buchberger algorithm provides a new set of rewrite rules $R'$, constructed from $R$ and a termination order, such that
\begin{enumerate}
\item $R'$ induces a confluent rewriting system,
\item the equivalence relations induced by $R'$ and $R$ are equal.
\end{enumerate}
Here, what we want to complete is a set $F$ of reduction operators. A \emph{completion of F} is a set $F'$ containing $F$ such that
\begin{enumerate}
\item $F'$ is confluent,
\item the two operators $\wedge F'$ and $\wedge F$ are equal.
\end{enumerate}
We show that a completion always exists. For that, we use the lattice structure to define an operator $C^F$ called the \emph{F-complement}. Our main result is Theorem~\ref{The F-completion is a completion} which states that the set $F\cup\left\{C^F\right\}$ is a completion of $F$. When $F$ is associated to a set of rewrite rules on $\K{X^*}$, the operator $C^F$ maps every obstruction $w$ to $\left(\wedge F\right)(w)$. In Theorem~\ref{Buchberger algorithm}, we use this operator to construct \G\ bases with reduction operators.

\paragraph{Reduction Operators without Total Order.}

The Knuth-Bendix completion algorithm does not require a total order on terms, which implies that it could fail. Indeed, at some point of the algorithm, one could have two normal forms $t_1$ and $t_2$ of a given term that we cannot compare for a fixed non-total order. The same phenomena holds for reduction operators when we do not assume that the order on $G$ is total. In this case, the restriction of the kernel map to reduction operators is not onto. In Section~\ref{Generalised reduction operators}, we deduce two important consequences of this fact. The first one is that it could happen that the lattice structure does not exist. The second one is more subtle: even if a set admits a lower bound, the latter does not necessarily have the "right shape". By right shape, we mean that this lower bound does not necessarily come from the lattice structure on the set of subspaces of $\K{G}$. As a consequence, the $F$-complement is not always defined. However, the existence of a lower bound with the right shape is sufficient to guarantee that it exists. We point out that Section~\ref{Generalised reduction operators} does not contain results non appearing in the previous sections, but is there to show that requiring a total order is of crucial importance in our constructions.

\paragraph{Reduction Operators and Computations.}

We explain the computational aspects of reduction operators. In order to do computations, we consider reduction operators relative to finite sets. This is probably a severe limitation. However, we explain in the conclusion of the paper why it should be possible, despite this limitation, to do effective computations. Fix a totally ordered finite set $\left(G,<\right)$. Since the lattice structure comes from the bijection between $\RO$ and subspaces of $\K{G}$, we need to have an explicit description of this bijection to compute the various algebraic constructions mentioned above. An online version of this implementation is available~\footnote{http://pastebin.com/0YZCfAD4}. This implementation uses the SageMath software~\footnote{http://www.sagemath.org}, written in Python. In the paper, we illustrate several constructions with examples. The examples for which we do not detail the computations were treated with the online version. Moreover, recall that reduction operators have applications to Koszulness. When an algebra has the Koszulness property, a family of important invariants, called \emph{homology groups} \cite{MR1438546}, can be expressed in terms of upper bound of reduction operators relative to finite sets. Hence, this is a context where our implementation can be used. The computation of these invariants already exists \cite{MR2030264}. Here, we propose an implementation using other techniques.

\begin{center}

\textbf{Organisation.}

\end{center}
In Section~\ref{Definition and lattice structure}, we define the notion of reduction operator relative to a well-ordered set. We also equip the set of these operators with a lattice structure and formulate the notion of confluence. In Section~\ref{F-compositions}, we define several notions from abstract rewriting theory in terms of reduction operators: normal forms, Church-Rosser property and local confluence. We show that the notions of confluence, local confluence and Church-Rosser property are equivalent. In Section~\ref{reduction operators and abstract rewriting}, we explain our notion of confluence from the viewpoint of abstract rewriting. Section~\ref{Confluence for a pair of reduction operators title} contains results about a pair of reduction operators. These results are necessary in Section~\ref{Completion}. In the latter, we define the notions of completion, complement and minimal complement. We also show that a minimal complement always exists. In Section~\ref{Presentations by operator}, we formulate the notions of presentation by operator and of confluent presentation by operator. We also link the latter with \G\ bases. In Section~\ref{Algebraic structure}, we formulate a general definition of reduction operator, which is relative to an ordered set. We show that the set of these operators is an ordered set but does not necessarily admit a lattice structure. In Section~\ref{Rewriting properties, general case}, we define the notion of completable set of reduction operators and study their rewriting properties. \newline\newline
\textbf{Acknowledgement.} This work was supported by the Sorbonne-Paris-Cit\'e IDEX grant Focal and the ANR grant ANR-13-BS02-0005-02 CATHRE. 

\section{Rewriting Properties of Reduction Operators}\label{Rewriting properties of reduction operators}

\subsection{Lattice Structure and Confluence}\label{Definition and lattice structure}

\subsubsection{Notations}\label{Notations}

We denote by $\K{}$ a commutative field. We say vector space instead of $\K{}$-vector space. Let $X$ be a set. We denote by $\K{X}$ the vector space with basis $X$: its non-zero elements are the finite formal linear combinations of elements of $X$ with coefficients in $\K{}$. An element of $X$ is called a \emph{generator} of $\K{X}$. By construction of $\K{X}$, for every $v\in\K{X}\setminus\{0\}$, there exist a unique finite subset $S_v$ of $X$ and a unique family of non zero scalars $\left(\lambda_x\right)_{x\in S_v}$ such that $v$ is equal to $\sum_{x\in S_v}\lambda_xx$. The set $S_v$ is the \emph{support} of $v$. 

\subsubsection{Leading Generator and Leading Coefficient}

Let $\left(G,<\right)$ be a \emph{well-ordered set}, that is, $G$ is a set and $<$ is a well-founded total order on $G$. The order on $G$ being total, every non-empty finite subset of $G$ admits a greatest element. In particular, for every $v\in\K{G}\setminus\{0\}$, the support of $v$ admits a greatest element, written $\lgen{v}$. We also write $\lc{v}=\lambda_{\lgen{v}}$. The elements $\lgen{v}$ and $\lc{v}$ are the \emph{leading generator} and the \emph{leading coefficient} of $v$, respectively. We extend the order $<$ on $G$ into a partial order on $\K{G}$ in the following way: we have $u<v$ if $u=0$ and $v$ is different from 0 or if $\text{lg}(u)<\text{lg}(v)$. 

Throughout Section~\ref{Rewriting properties of reduction operators} we fix a well-ordered set $\left(G,<\right)$.

\subsubsection{Reduction Operators}\label{Definition of reduction operators}

A \emph{reduction operator relative to} $\left(G,<\right)$ is an idempotent linear endomorphism $T$ of $\K{G}$ such that for every $g\in G$, we have $T(g)\leq g$. We denote by $\textbf{RO}\left(G,<\right)$ the set of reduction operators relative to $\left(G,<\right)$. Given $T\in\RO$, a generator $g$ is said to be \emph{T-reduced} if $T(g)$ is equal to $g$. We denote by $\red{T}$ the set of $T$-reduced generators and by $\nred{T}$ the complement of $\red{T}$ in $G$.

\subsubsection{Remarks}
Let $T\in\RO$.
\begin{enumerate}
\item The image of $T$ is the vector space spanned by $T$-reduced generators:
\[\im{T}=\K{\red{T}}.\]
\item Let $g\in G$. The condition $T(g)\leq g$ means that one of the following two conditions is fulfilled:
\begin{enumerate}
\item $g$ is $T$-reduced,
\item $T(g)$ is a linear combination of elements of $G$ strictly smaller than $g$ for $<$.
\end{enumerate}
\end{enumerate}

\subsubsection{Reduction Matrices}\label{Reduction matrix}

In our examples, we sometimes consider the case where $\left(G,<\right)$ is a totally ordered finite set: $G=\{g_1<\cdots<g_n\}$. In this case, we use matrix notations to describe linear maps, and thus, to describe reduction operators. For that, given an endomorphism $T$ of $\K{G}$, the matrix of $T$ with respect to the basis $\{g_1,\cdots ,g_n\}$ is called the \emph{canonical matrix of T relative to $\left(G,<\right)$.} We consider the convention that the $j$-th column of this matrix contains the coefficients of $T(g_j)$ with respect to the basis $G$. Moreover, we say that a square matrix $M$ is a \emph{reduction matrix} if the following conditions are fulfilled:
\begin{enumerate}
\item\label{Upper triangular} $M$ is upper triangular and the elements of its diagonal are equal to 0 or 1,
\item\label{Matrice de red non-red} if an element of the diagonal of $M$ is equal to 0, then the other elements of the line to which it belongs are equal to 0,
\item\label{Matrice de red réduit} if an element of the diagonal of $M$ is equal to 1, then the other elements of the column to which it belongs are equal to 0.
\end{enumerate}

Our purpose is to show that an endomorphism of $\K{G}$ is a reduction operator relative to $\left(G,<\right)$ if and only if its canonical matrix relative to $\left(G,<\right)$ is a reduction matrix. For that we need the following:

\subsubsection{Lemma}\label{A reduction matrix is idempotent}

\emph{A reduction matrix is idempotent.}

\begin{proof}

Let $M$ be a reduction matrix. Let $\left(m_{ij}\right)_{1\leq i,j\leq n}$ be the coefficients of $M$, where $i$ and $j$ denote the row and the column of $m_{ij}$, respectively. Let $A=\left(a_{ij}\right)_{1\leq i,j\leq n}$ be the product $M\times M$. For every $1\leq i,j\leq n$, we have
\[a_{ij}=\sum_{k=1}^nm_{ik}m_{kj}.\]
Let $1\leq i\leq n$ such that $m_{ii}=0$. From Point~\ref{Matrice de red non-red} of~\ref{Reduction matrix}, for every $1\leq k\leq n$, we have $m_{ik}=0$. Thus, for every $1\leq j\leq n$, we have
\[\begin{split}
a_{ij}&=0\\
&=m_{ij}.
\end{split}\]
Hence, the $i$-th rows of $M$ and of $A$ are equal when $m_{ii}$ is equal to $0$. Let $1\leq i\leq n$ such that $m_{ii}=1$. For every $1\leq j\leq n$, we have
\[\begin{split}
a_{ij}&=\sum_{k\neq i}m_{ik}m_{kj}+m_{ii}m_{ij}\\
&=\sum_{k\neq i}m_{ik}m_{kj}+m_{ij}.
\end{split}\]
Let $k\neq i$ such that $m_{ik}$ is different form $0$. From Point~\ref{Matrice de red réduit} of~\ref{Reduction matrix}, $m_{kk}$ is different from $1$ so that it is equal to $0$. Thus, from Point~\ref{Matrice de red non-red}, $m_{kj}$ is equal to $0$. Thus, $m_{ik}m_{kj}$ is equal to $0$ for every $k\neq i$, so that $a_{ij}$ is equal to $m_{ij}$. Hence, the $i$-th rows of $M$ and of $A$ are equal when $m_{ii}$ is equal to $1$. Hence, the rows of $M$ and $A$ are equal so that $A$ and $M$ are equal, that is, $M$ is idempotent.

\end{proof}

\subsubsection{Proposition}

\emph{Assume that G is finite. An endomorphism of $\K{G}$ is a reduction operator relative to $\left(G,<\right)$ if and only if its canonical matrix relative to $\left(G,<\right)$ is a reduction matrix.}

\begin{proof}

Let $T$ be an endomorphism of $\K{G}$ and let $M$ be its canonical matrix relative to $\left(G,<\right)$.

Assume that $T$ belongs to $\RO$. For every $g\in G$, we have $T(g)\leq g$. In particular, $M$ satisfies Point~\ref{Upper triangular} and Point~\ref{Matrice de red réduit} of~\ref{Reduction matrix}.Moreover, the image of $T$ being equal to the vector space spanned by $\red{T}$, no element of $\nred{T}$ belongs to the decomposition of an element of $T(g)$ for $g\in G$. Hence, $M$ satisfies~\ref{Matrice de red non-red} of~\ref{Reduction matrix}. Thus, $M$ is a reduction matrix.

Assume that $M$ is a reduction matrix. From Point~\ref{Upper triangular} and Point~\ref{Matrice de red réduit} of~\ref{Reduction matrix}, for every $g\in G$, we have $T(g)\leq g$. Moreover, from Lemma~\ref{A reduction matrix is idempotent}, $M$ is idempotent so that $T$ is idempotent. Hence, $T$ is a reduction operator relative to $\left(G,<\right)$.

\end{proof}

\subsubsection{Examples}\label{Example of Reduction matrices}

\begin{enumerate}
\item The zero matrix and the identity matrix are reduction matrices. More generally, every diagonal matrix admitting only 0 or 1 on its diagonal is a reduction matrix.
\item The matrix
\[\begin{pmatrix}
1&1&1&1\\
0&0&0&0\\
0&0&0&0\\
0&0&0&0
\end{pmatrix},\]
is a reduction matrix.
\item\label{Definition of the braided monoid} From~\cite{MR785160}, we consider the presentation of the monoid ${\textbf{B}_3}^+$ by 3 generators $x,\ y$ and $z$ subject to two relations $yz=x$ and $zx=xy$. Let $X=\{x<y<z\}$. Let $X^*$ be the free monoid over $X$: this is the set of (possibly empty) words written with the alphabet $X$. This set is totally ordered for the deg-lex order, still denoted by $<$, induced by the order on $X$. For this order, $yz$ is greater than $x$ and $zx$ is greater than $xy$. Hence, the endomorphism of $\K{X^*}$ defined on the basis $X^*$ by $S(yz)=x$, $S(zx)=xy$ and $S(w)=w$ for every $w\in X^*\setminus\{yz,zx\}$, is a reduction operator relative to $\left(X^*,<\right)$.
\end{enumerate}

\subsubsection{Reduced Basis}\label{Normalised echelon basis}

Let $V$ be a subspace of $\K{G}$. A \emph{reduced basis} of $V$ is a basis $\M{B}$ of $V$ satisfying the following conditions:
\begin{enumerate}
\item\label{Normalised condition} for every $e\in\M{B}$, $\lc{e}$ is equal to 1,
\item\label{Echelon condition} given two different elements $e$ and $e'$ of $\M{B}$, $\lgen{e'}$ does not belong to the support of $e$.
\end{enumerate}

\subsubsection{Notation}\label{Notation of reduced basis}

Let $V$ be a subspace of $\K{G}$ and let $\M{B}$ be a reduced basis of $V$. Let $e$ and $e'$ be two distinct elements of $\M{B}$. The condition~\ref{Echelon condition} of~\ref{Normalised echelon basis} implies that $\lgen{e}$ is different from $\lgen{e'}$. Hence, $\M{B}$ is indexed by the set $\tilde{G}=\left\{\lgen{e}\ \mid\ e\in\M{B}\right\}$. In the sequel, a reduced basis $\M{B}$ is written $\M{B}=\left(e_g\right)_{g\in\tilde{G}}$, where for every $g\in\tilde{G}$, we have $\lgen{e_g}=g$.

\subsubsection{Remark}\label{Remark on n.e.b}

Let $V$ be a subspace of $\K{G}$ and let $\left(e_g\right)_{g\in\tilde{G}}$ be a reduced basis of $V$. For every $g\in G$, let $V_g$ be the set of elements of $V$ with leading generator $g$. The set $V_g$ is non empty if and only if $g$ belongs to $\tilde{G}$. Hence, if $\M{B}_1=\left(e_g\right)_{g\in\tilde{G}_1}$ and $\M{B}_2=\left(e'_g\right)_{g\in\tilde{G}_2}$ are two reduced bases of $V$, the two sets $\tilde{G}_1$ ant $\tilde{G}_2$ are equal.

\subsubsection{Examples}\label{Example of normalised echelon basis}

\begin{enumerate}
\item\label{Computation of a reduced basis} Let $G=\{g_1<g_2<g_3<g_4\}$ and let $V$ be the subspace of $\K{G}$ spanned by the elements $v_1=g_2-g_1$, $v_2=g_4-g_3$ and $v_3=g_4-g_2$. The elements $\lgen{v_2}$ and $\lgen{v_3}$ are equal to $g_4$, so that $\{v_1,v_2,v_3\}$ is not a reduced basis of $V$. Letting ${v_2}'=v_3-v_2$, that is, ${v_2}'=g_3-g_2$, the set $\{v_1,{v_2}',{v_3}\}$ is a basis of $V$. This is still not a reduced basis because $\lgen{v_1}=g_2$ appears in the supports of ${v_2}'$ and ${v_3}$. Letting ${v_2}''={v_2}'+v_1$ and ${v_3}'=v_3+v_1$, that is, ${v_2}''=g_3-g_1$ and ${v_3}'=g_4-g_1$, the set $\{v_1,{v_2}'',{v_3}'\}$ is a basis of $V$. This basis is reduced. Using the notation introduced in~\ref{Notation of reduced basis}, it is equal to $\left\{e_{g_2},e_{g_3},e_{g_4}\right\}$, where $e_{g_i}$ is equal to $g_i-g_1$ for $i\in\{2,3,4\}$.
\item A \emph{semi-reduced} string rewriting system is a string rewriting system $\EV{X\mid R}$ such that the left-hand sides of the elements of $R$ are pairwise distinct and if no right-hand side of an element of $R$ is the left-hand side of another one. The set $\{w-w'\ \mid\ w\F{}w'\in R\}$ is a reduced basis of the vector space it spans.
\end{enumerate}

\subsubsection{Theorem}\label{Unicity of normalised echelon basis}

\emph{Let $\left(G,<\right)$ be a well-ordered set. Every subspace of $\K{G}$ admits a unique reduced basis.}

\begin{proof}

Let $V$ be a subspace of $\K{G}$. First, we construct by induction on $G$ a reduced basis of $V$. Let $g_0$ be the smallest element of $G$. If $V_{g_0}$ is empty, we let $\M{B}_{g_0}=\emptyset$. In the other case, $g_0$ belongs to $V$ and we let $\M{B}_{g_0}=\{g_0\}$. Let $g\in G$. Assume by induction that for every $g'<g$ we have built a set $\M{B}_{g'}$ such that the following conditions are fulfilled:
\begin{enumerate}\label{At most one element}
\item For every $g'<g$, the set $\M{B}_{g'}$ contains at most one element. 
\item\label{Second point of induction} Let
\[I_g=\left\{g'\in G\ \mid\ g'<g\ \text{and}\ \M{B}_{g'}\neq\emptyset\right\}.\]
For every $g'\in I_g$
\begin{enumerate}
\item\label{1} the unique element $e_{g'}$ of $\M{B}_{g'}$ belongs to $V$,
\item\label{2} $\lgen{e_{g'}}$ is equal to $g'$ and $\lc{e_{g'}}$ is equal to $1$,
\item\label{3} for every $\tilde{g}\in I_g$ such that $\tilde{g}$ is different from $g'$, ${\tilde{g}}$ does not belong to the support of $e_{g'}$,
\item\label{4} the set $V_{g'}$ is included in $\K{\left\{e_{\tilde{g}}\ \mid\ \tilde{g}\in I_g\right\}}$.
\end{enumerate}
\end{enumerate}
If $V_g$ is empty, we let $\M{B}_{g}=\emptyset$. If $V_{g}$ is non empty, let $v_{g}$ be an element of $V_{g}$ such that $\lc{v_{g}}$ is equal to 1. In particular, $v_g$ admits a decomposition
\[v_{g}=\sum_{g'\in J}\mu_{g'}g'+g,\]
where for every $g'\in J$, we have $g'<g$. We let $\M{B}_{g}=\{e_{g}\}$, where
\[e_{g}=v_{g}-\sum_{g'\in I_g}\mu_{g'}e_{g'}.\]
In particular, $\M{B}_g$ contains at most one element. By construction, $e_{g}$ belongs to $V$, $\lgen{e_{g}}$ is equal to $g$ and $\lc{e_{g}}$ is equal to 1, so that Point~\ref{1} and Point~\ref{2} hold. Moreover, for every $g'\in I_g$, $g'$ does not belong to the support of $e_{g}$, so that Point~\ref{3} holds. It remains to show that $V_g$ is included in the vector space spanned by the elements $e_{g'}$ where $g'$ belongs to $I_g\cup\{g\}$. This vector space is equal to $\K{\left\{e_{\tilde{g}}\ \mid\ \tilde{g}\in I_g\right\}}\oplus\K{e_g}$. Let $v$ be an element of $V_g$. Then, $v-\lc{v}e_{g}$ belongs to $V$ and $\lgen{v-\lc{v}e_{g}}$ is strictly smaller than $g$. By the induction hypothesis, $v-\lc{v}e_{g}$ belongs to $\K{\left\{e_{\tilde{g}}\ \mid\ \tilde{g}\in I_g\right\}}$, so that $v$ belongs to $\K{\left\{e_{\tilde{g}}\ \mid\ \tilde{g}\in I_g\right\}}\oplus\K{e_g}$. This inductive construction provides a family of sets $\left(\M{B}_g\right)_{g\in G}$ such that $\M{B}=\bigcup_{g\in G}\M{B}_g$ is a generating set of $V=\bigsqcup_{g\in G}V_g$. Moreover, this family is free because the leading generators of its elements are pairwise distinct. Hence, $\M{B}$ is a basis of $V$. This basis is reduced by construction.

Let us show that such a basis is unique. Let $\M{B}_1=\left(e_g\right)_{g\in\tilde{G}_1}$ and $\M{B}_2=\left(e'_g\right)_{\in\tilde{G}_1}$ be two reduced bases of $V$. We have seen in Remark~\ref{Remark on n.e.b} that $\tilde{G}_1$ and $\tilde{G}_2$ are equal, so that we write $\M{B}_1=\left(e_g\right)_{g\in\tilde{G}}$ and $\M{B}_2=\left(e'_g\right)_{\in\tilde{G}}$. Let $g\in\tilde{G}$. Assume that $e_g$ is different from $e'_g$. Let
\[e'_g-e_g=\sum_{g'\in I}\lambda_{g'}e'_{g'},\]
be the decomposition of $e'_g-e_g$ with respect to the basis $\M{B}_2$. The leading generator of ${e'_g-e_g}$ is equal to the greatest element of $I$, so that it belongs to $\tilde{G}$. Let $\tilde{G}^c$ be the complement of $\tilde{G}$ in $G$. The condition~\ref{Echelon condition} of the definition of a reduced basis implies that $e_g-g$ and $e'_g-g$ belong to $\K{\tilde{G}^c}$. Hence, $\lgen{e'_g-e_g}$ belongs to $\tilde{G}^c$, which is contradiction. Thus, for every $g\in\tilde{G}$, the two elements $e_g$ and $e'_g$ are equal, that is, $\M{B}_1$ and $\M{B}_2$ are equal.

\end{proof}

\subsubsection{Proposition}\label{Bijection}

\emph{Let} $V$ \emph{be a subspace of $\K{G}$. There exists a unique reduction operator T with kernel V. Moreover,} $\nred{T}$ \emph{is equal to $\tilde{G}$, where $\left(e_g\right)_{g\in\tilde{G}}$ is the reduced basis of} $V$.

\begin{proof}

Let $T$ be the endomorphism of $\K{G}$ defined on the basis $G$ in the following way
\[T(g)=\left\{
\begin{split}
& g-e_g,\ \text{if}\ g\in\tilde{G},\\
&g,\ \text{ohterwise}.
\end{split}\right.\]
By definition of a reduced basis, $T$ is a reduction operator relative to $\left(G,<\right)$. By construction, the kernel of $T$ is equal to $V$ and $\nred{T}$ is equal to $\tilde{G}$.

Let $T_1$ and $T_2$ be two reduction operators with kernel $V$. The two sets 
\[\{g-T_1(g)\ \mid\ T_1(g)\neq g\}\ \ \text{and}\ \ \{g-T_2(g)\ \mid\ T_2(g)\neq g\},\]
are reduced bases of $V$. From Theorem~\ref{Unicity of normalised echelon basis}, these sets are equal. Hence, $T_1(g)$ is different from $g$ if and only if $T_2(g)$ is different from $g$ and in this case, $T_1(g)$ is equal to $T_2(g)$. It follows that $T_1(g)$ and $T_2(g)$ are equal for every $g\in G$, so that $T_1$ and $T_2$ are equal.
\end{proof}

\subsubsection{Notation}

Proposition~\ref{Bijection} implies that the kernel map induces a bijection between $\RO$ and the set of subspaces of $\K{G}$. The inverse of this bijection is written $\theta$.

\subsubsection{Lattice Structure}\label{Lattice structure}

We consider the binary relation on $\RO$ defined by
\[T_1\preceq T_2\ \text{if and only if}\ \ker\left(T_2\right)\ \subseteq\ \ker\left(T_1\right).\]
This relation is reflexive and transitive. From Proposition~\ref{Bijection}, it is also anti-symmetric. Hence, it is an order relation on $\RO$. Moreover, we have the equivalence:
\begin{equation}\label{Characterisation of the order relation}
T_1\preceq T_2\ \text{if and only if}\ T_1\circ T_2=T_1.
\end{equation}
Let us equip $\RO$ with a lattice structure. The lower bound $T_1\wedge T_2$ and the upper bound $T_1\vee T_2$ of two elements $T_1$ and $T_2$ of $\RO$ are defined in the following manner:
\[\left\{\begin{split}
&T_1\wedge T_2=\theta\left(\ker(T_1)+\ker(T_2)\right), \\ 
&T_1\vee T_2=\theta\left(\ker(T_1)\cap\ker(T_2)\right).
\end{split}
\right.\]

\subsubsection{Example}\label{Illustration of the lattice structure}

Let $G=\left\{g_1<g_2<g_3<g_4\right\}$ and $P=\left(T_1,T_2\right)$, where
\[T_1=\begin{pmatrix}
1&1&0&0\\
0&0&0&0\\
0&0&1&1\\
0&0&0&0
\end{pmatrix}\ \text{and}\ 
T_2=\begin{pmatrix}
1&0&0&0\\
0&1&0&1\\
0&0&1&0\\
0&0&0&0
\end{pmatrix}.\]
The kernels of $T_1$ and $T_2$ are equal to $\K{\{g_2-g_1\}}\oplus\K{\{g_4-g_3\}}$ and $\K{\{g_4-g_2\}}$, respectively. The kernel of $T_1\wedge T_2$ is the vector space spanned by $v_1=g_2-g_1$, $v_2=g_4-g_3$ and $v_3=g_4-g_2$. This is the vector space of Example~\ref{Example of normalised echelon basis}, Point~\ref{Computation of a reduced basis}. Hence, the kernel of $T_1\wedge T_2$ is the vector space spanned by $\{g_2-g_1,g_3-g_1,g_4-g_1\}$, so that we have
\[T_1\wedge T_2=\begin{pmatrix}
1&1&1&1\\
0&0&0&0\\
0&0&0&0\\
0&0&0&0
\end{pmatrix}.\]

\subsubsection{Lemma}\label{Inclusion of images}

\emph{Let $T_1$ and $T_2$ be two reduction operators relative to $\left(G,<\right)$. Then, we have:}
\[T_1\preceq T_2\Longrightarrow\red{T_1}\ \subseteq\ \red{T_2}.\]

\begin{proof}

Let $\M{B}_1=\left(e_g\right)_{g\in\tilde{G}_1}$ and $\M{B}_2=\left(e'_g\right)_{g\in\tilde{G}_2}$ be the reduced bases of $\ker\left(T_1\right)$ and $\ker\left(T_2\right)$, respectively. We consider the notations of Remark~\ref{Remark on n.e.b}: given a subspace $V$ of $\K{G}$, let $V_g$ be the set of elements of $V$ with leading generator $g$. The sets
\[\left\{g\in G\ \mid\ \ker\left(T_1\right)_g\neq\emptyset\right\}\ \ \text{and}\ \ \left\{g\in G\ \mid\ \ker\left(T_2\right)_g\neq\emptyset\right\},\]
are equal to $\tilde{G}_1$ and $\tilde{G}_2$, respectively. Hence, if $\ker\left(T_2\right)$ is included in $\ker\left(T_1\right)$, then $\tilde{G}_2$ is included in $\tilde{G}_1$. From Proposition~\ref{Bijection}, we deduce that $\red{T_1}$ is included in $\red{T_2}$.

\end{proof}

\subsubsection{Obstructions}\label{Obstructions}

Let $F$ be a subset of $\RO$. We let
\[\red{F}=\bigcap_{T\in F}\red{T}\ \text{and}\ \wedge F=\theta\left(\sum_{T\in F}\ker\left(T\right)\right).\]
For every $T\in F$, we have $\wedge F\preceq T$. Thus, from Lemma~\ref{Inclusion of images}, the set $\red{\wedge F}$ is included in $\red{T}$ for every $T\in F$, so that it is included in $\red{F}$. We write
\begin{equation} \label{Definition of obstructions}
\obsred{F}=\red{F}\setminus\red{\wedge F}.
\end{equation}

\subsubsection{Confluence}\label{Confluence of a finite family of reduction operators}

A subset $F$ of $\RO$ is said to be \emph{confluent} if $\obsred{F}$ is the empty set.

\subsubsection{Examples}\label{Examples for confluence}

\begin{enumerate}
\item\label{Example of confluence, obstructions} We consider the pair $P$ of Example~\ref{Illustration of the lattice structure}. The set $\red{P}$ is equal to $\{g_1,g_3\}$. Moreover, $\red{T_1\wedge T_2}$ is equal to $\{g_1\}$, so that $\obsred{P}$ is equal to $\{g_3\}$. Hence, $P$ is not confluent.
\item\label{Non confluence for the braided monoid} We consider the operator $S$ defined in Point~\ref{Definition of the braided monoid} of Example~\ref{Example of Reduction matrices}. Let $S_1$ and $S_2$ be the restrictions to the vector space spanned by words of length smaller or equal to $3$ of $S\otimes\id{\K{X}}$ and $\id{\K{X}}\otimes S$, respectively. These two operators are defined for every word $w$ of length smaller or equal to $3$ by
\begin{enumerate}
\item $S_1(yzt)=xt$, $S_1(zxt)=xyt$ for every $t\in X$ and $S_1(w)=w$ is different from $yzt$ and $zxt$ for some $t\in X$,
\item $S_2(tyz)=tx$, $S_2(tzx)=zxy$ for every $t\in X$ and $S_2(w)=w$ is different from $tyz$ and $zzx$ for some $t\in X$.
\end{enumerate}
Let $P=\left(S_1,S_2\right)$. The set of not $\wedge P$-reduced generators is equal to \[\Big\{xyz,xzx,yxy,yyz,yzx,yzy,yzz,zxx,zxy,zxz,zyz,zzx\Big\}.\] 
We remark that $yxy$ belongs to this list but also belongs to $\red{P}$ since its two sub-words of length $2$ are $S$-reduced. Hence, $yxy$ belongs to $\obsred{P}=\red{P}\setminus\red{\wedge P}$, so that $P$ is not confluent.
\end{enumerate}

\subsection{Normal Forms, Church-Rosser Property and Newman's Lemma}\label{F-compositions}

Throughout this section we fix a subset $F$ of $\RO$. We denote by $\EV{F}$ the submonoid of $\left(\text{End}\left(\K{G}\right),\circ\right)$ spanned by $F$. 

\subsubsection{Multi-Set Order}

Given an element $v$ of $\K{G}$, let $S_v$ be the support of $v$. We introduce the order $\leq_\text{mul}$ on $\K{G}$ defined in the following way: we have $v\leq_\text{mul}v'$ if for every $g\in S_v$ such that $g$ does not belong to $S_{v'}$, there exists an element $g'\in S_{v'}$ not appearing in $S_v$, such that $g<g'$. The order $<$ being well-founded, this is also the case for $\leq_\text{mul}$. Moreover, given $v\in\K{G}$ and $T\in F$ such that $v$ does not belong to $\K{\red{T}}$, $T(v)$ is strictly smaller than $v$ for $\leq_\text{mul}$. Finally, we also remark that 0 is the smallest element of $\K{G}$ for $\leq_\text{mul}$.

\subsubsection{Normal Forms}

\begin{enumerate}
\item An \emph{F-normal form} is an element of $\K{\red{F}}$.
\item Let $v$ and $v'$ be two elements of $\K{G}$. We say that \emph{v rewrites into} $v'$ if there exists $R\in\EV{F}$ such that $v'$ is equal to $R(v)$.
\item Let $v$ be an element of $\K{G}$. An \emph{F-normal form of v} is an $F$-normal form $v'$ such that $v$ rewrites into $v'$.
\end{enumerate}

\subsubsection{Remark}

Let $v$ and $v'$ be two elements of $\K{G}$ such that $v$ rewrites into $v'$. Reasoning by induction on the length of $R\in\EV{F}$ such that $v'=R(v)$, we conclude that $v-v'$ belongs to $\ker\left(\wedge F\right)$.

\subsubsection{Church-Rosser Property}

We say that $F$ has the \emph{Church-Rosser property} if for every $v\in\K{G}$, $v$ rewrites into $\left(\wedge F\right)(v)$.

\subsubsection{Theorem}\label{Théorème de Church-Rosser}

\emph{A subset of $\emph{\textbf{RO}}\left(G,<\right)$ is confluent if and only if it has the Church-Rosser property.}

\begin{proof}

Let $F$ be a subset of $\RO$.

Assume that $F$ has the Church-Rosser property. Let $g$ be an element of $\red{F}$. For every $T\in F$, $T(g)$ is equal to $g$. As a consequence, $R(g)$ is equal to $g$ for every $R\in\EV{F}$. The set $F$ having the Church-Rosser property, we deduce that $g$ is equal to $\left(\wedge F\right)(g)$, so that $g$ belongs to $\red{\wedge F}$. That shows that $\red{F}$ is included in $\red{\wedge F}$, that is, $F$ is confluent.

Assume that $F$ is confluent. Let us show by induction on $\leq_\text{mul}$ that for every $v\in\K{G}$, $v$ rewrites into $\left(\wedge F\right)(v)$. If $v$ is equal to $0$, this is obvious. Let $v\in\K{G}$. Assume by induction that for every $v'\in\K{G}$ such that $v'\leq_\text{mul}v$, $v'$ rewrites into $\left(\wedge F\right)(v')$. If $v$ belongs to $\K{\red{\wedge F}}$, then $v$ is equal to $\left(\wedge F\right)(v)$, so that $v$ rewrites into $\left(\wedge F\right)(v)$. Assume that $v$ does not belong to $\K{\red{\wedge F}}$. The set $F$ being confluent, $v$ does not belong to $\K{\red{F}}$, that is, there exists $T\in F$ such that $T(v)$ is different from $v$. The element $T(v)$ is strictly smaller than $v$ for $\leq_\text{mul}$. By the induction hypothesis, there exists $R\in\EV{F}$ such that $R(T(v))$ is equal to $\left(\wedge F\right)\left(T(v)\right)$. The inequality $\wedge F\preceq T$ implies that $\wedge F\circ T$ is equal to $\wedge F$, so that $\left(\wedge F\right)\left(T(v)\right)$ is equal to $\left(\wedge F\right)(v)$. Hence, $R'=R\circ T$ is an element of $\EV{F}$ such that $R'(v)$ is equal to $\left(\wedge F\right)(v)$. Thus, $v$ rewrites into $\left(\wedge F\right)(v)$.

\end{proof}

\subsubsection{Lemma}\label{Terminaison}

\emph{Let v be an element of $\K{G}$ and let $(R_1,\cdots,R_n,\cdots)$ be a sequence of elements of $\EV{F}$ such that for every integer n, $R_{n}$ is a right divisor of $R_{n+1}$ in $\EV{F}$. The sequence $\left(R_n(v)\right)_{n\in\N}$ is stationary.}

\begin{proof}

We proceed by induction on $v$. If $v$ is equal to 0, then the sequence $\left(R_n(v)\right)_{n\in\N}$ is constant, equals to 0. Let $v\in\K{G}$. Assume that Lemma~\ref{Terminaison} holds for every $v'\in\K{G}$ such that $v'$ is strictly smaller than $v$ for $\leq_\text{mul}$. If the sequence $\left(R_n(v)\right)_{n\in\N}$ is constant equals to $v$, there is nothing to prove. In the other case, there exists $n_0$ such that $R_{n_0}(v)$ is different from $v$, so that we have $R_{n_0}(v)\leq_\text{mul}v$. By hypothesis, for every integer $n$, $R_{n}$ is a right divisor of $R_{n+1}$ in $\EV{F}$, that is, there exists $R'_n\in\EV{F}$ such that $R_{n+1}$ is equal to $R'_n\circ R_{n}$. Let $\left(Q_n\right)_{n\in\N}$ be the sequence of elements of $\EV{F}$ defined by
\[Q_1=R'_{n_0}\ \text{and}\ Q_{n+1}=R'_{n_0+n}\circ Q_n.\]
For every integer $n$, $Q_{n}$ is a right divisor of $Q_{n+1}$ in $\EV{F}$. By the induction hypothesis, the sequence $\left(Q_n\left(R_{n_0}(v)\right)\right)_{n\in\N}$ is stationary. Moreover, for every integer $n$, $Q_n\circ R_{n_0}$ is equal to $R_{n_0+n}$, so that the sequence $\left(R_n(v)\right)_{n\in\N}$ is stationary. Hence, Lemma~\ref{Terminaison} holds.

\end{proof}

\subsubsection{Proposition}\label{Existence d'une forme normale}

\emph{Every element of $\K{G}$ admits an F-normal form.}

\begin{proof}

Let $v$ be an element of $\K{G}$. We have to show that there exists $R\in\EV{F}$ such that $R(v)$ belongs to $\K{\red{F}}$. Assume by way of contradiction that for every $R\in\EV{F}$, $R(v)$ does not belong to $\K{\red{F}}$. The morphism $\id{\K{G}}$ belonging to $\EV{F}$, $v$ does not belong to $\K{\red{F}}$. In particular, there exists $T_1\in F$ such that $v$ does not belong to $\K{\red{T_1}}$. Assume that we have constructed elements $T_1,\cdots,T_n$ of $F$. The morphism $R_n=T_n\circ\cdots\circ T_1$ belongs to $\EV{F}$. Hence, the element $R_n(v)$ does not belong to $\K{\red{F}}$, so that there exists $T_{n+1}~\in~ F$ such that $R_n(v)$ does not belong to $\K{\red{T_{n+1}}}$. This process enables us to obtain a sequence $\left(R_n\right)_{n\in\N}$ of elements of $\EV{F}$ such that for every integer $n$, $R_{n}$ is a right divisor of $R_{n+1}$ in $\EV{F}$ and such that the sequence $\left(R_n(v)\right)_{n\in\N}$ is not stationary. This is a contradiction with Lemma~\ref{Terminaison}. Hence, Proposition~\ref{Existence d'une forme normale} holds.

\end{proof}

\subsubsection{Notation}\label{Minimal element}

For every $v\in\K{G}$, let $[v]$ be the set of elements $v'\in\K{G}$ such that $v'-v$ belongs to $\ker\left(\wedge F\right)$ . 

\subsubsection{Lemma}\label{La borne inf est l'opérateur de normalisation}

\emph{For every $v\in\K{G}$, $\left(\wedge F\right)(v)$ is the smallest element of $[v]$ for the order $<_\emph{mul}$. Moreover, if every element v of $\K{G}$ admits exactly one F-normal form, this normal form is equal to $\left(\wedge F\right)(v)$. }

\begin{proof}

Let us show the first assertion. Let $v\in\K{G}$ and let $v'\in [v]$. In particular, $v'$ belongs to $[\left(\wedge F\right)(v)]$, that is, there exists $v''\in\ker\left(\wedge F\right)$ such that
\begin{equation}\label{Decomposition of v'}
v'=\left(\wedge F\right)(v)+v''.
\end{equation}
The element $v''$ belonging to $\ker\left(\wedge F\right)$, it admits a decomposition 
\begin{equation}\label{Decomposition of v''}
v''=\sum\lambda_i\left(g_i-\left(wedge F\right)(g_i)\right),
\end{equation}
where each $g_i$ is not $\wedge F$-reduced. Let $g\in G$ be an element of the support of $\left(\wedge F\right)(v)$ not appearing in the one of $v'$. Let us show that there exists an index $i$ such that
\begin{enumerate}
\item $g_i$ is strictly greater than $g$,
\item $g_i$ does not belong to the support of $\left(\wedge F\right)(v)$,
\item $g_i$ belongs to the support of $v'$.
\end{enumerate}
Relation (\ref{Decomposition of v'}) and the hypothesis on $g$ imply that the latter belongs to the support of $v''$. Moreover, $g$ belongs to the image of $\wedge F$, that is, it is $\wedge F$-reduced. From Relation (\ref{Decomposition of v''}), we deduce that $g$ belongs to the support of $\left(\wedge F\right)(g_i)$ for some $i$. The element $g_i$ being not $\wedge F$-reduced, $g$ is strictly smaller than $g_i$ and does not belong to the support of $\left(\wedge F\right)(v)$. Finally, $g_i$ belonging to the support of $v''$ and not to the one of $\left(\wedge F\right)(v)$, Relation (\ref{Decomposition of v'}) implies that it belongs to the support of $v'$. Hence, we have $\left(\wedge F\right)(v)\leq_\text{mul}v'$ for every $v'\in[v]$, so that the first part of the lemma holds.

Let us show the second part of the lemma. Assume that every element $v$ of $\K{G}$ admits a unique $F$-normal form, written $N(v)$. It is clear that the operator $N$ is idempotent. Moreover, for every $R\in\EV{F}$ and for every $g\in G$, we have either $R(g)=g$ or $R(g)<g$. As a consequence, for every $g\in G$, we have either $N(g)=g$ or $N(g)<g$. We conclude that $N$ belongs to $\RO$. Let us show that $N$ is equal to $\wedge F$. Let $T$ be an element of $F$ and let $v$ be an element of $\ker\left( T\right)$. The element $v$ rewrites into $0$, so that $N(v)$ is equal to 0. Thus, the kernel of $T$ is included in the kernel of $N$ for every $T\in F$, that is, $N$ is smaller or equal to $T$ for every $T\in F$. Thus, we have the inequality $N\preceq\wedge F$. Moreover, from Relation~(\ref{Characterisation of the order relation}) (see~\ref{Lattice structure}), for every $T\in F$, the operator $\wedge F\circ T$ is equal to $\wedge F$. Hence, for every $R\in\EV{F}$, the operator $\wedge F\circ R$ is equal to $\wedge F$. As a consequence, for every $v\in\K{G}$, $\left(\wedge F\circ N\right)(v)$ being equal to $\left(\wedge F\circ R\right)(v)$ for some $R\in\EV{F}$, $\wedge F\circ N$ is equal to $\wedge F$. Using again Relation~(\ref{Characterisation of the order relation}), we have $\wedge F\preceq N$. Hence, $N$ is equal to $\wedge F$, so that Lemma~\ref{La borne inf est l'opérateur de normalisation} holds. 

\end{proof}

\subsubsection{Proposition}\label{Confluence vs unique forme normale}

\emph{The set F is confluent if and only if every element of $\K{G}$ admits a unique F-normal form.}

\begin{proof}

Assume that $F$ is confluent. Let $v$ be an element of $\K{G}$ and let $v_1$ and $v_2$ be two $F$-normal forms of $v$. Let $R_1$ and $R_2$ be two elements of $\EV{F}$ such that $R_i(v)$ is equal to $v_i$, for $i=1$ or 2. The elements $v-v_1$ and $v-v_2$ belong to $\ker\left(\wedge F\right)$. Hence, $v_1-v_2$ belongs to $\ker\left(\wedge F\right)$. Moreover, $v_1$ and $v_2$ belonging to $\K{\red{F}}$, $v_1-v_2$ also belongs to $\K{\red{F}}$, that is, $\K{\red{\wedge F}}$ since $F$ is confluent. Thus, $v_1-v_2$ belongs to the vector space $\K{\red{\wedge F}}~\cap~\ker\left(\wedge F\right)$. The operator $\wedge F$ being a projector, this vector space is reduced to $\{0\}$. We conclude that $v_1$ is equal to $v_2$, so that $v$ admits a unique $F$-normal form.

Assume that every element of $\K{G}$ admits a unique $F$-normal form. Let $v$ be an element of $\K{G}$. From Lemma~\ref{La borne inf est l'opérateur de normalisation}, the normal form of $v$ is equal to $\left(\wedge F\right)(v)$. Hence, $v$ rewrites into $\left(\wedge F\right)(v)$. We conclude that $F$ has the Church Rosser-property, that is, $F$ is confluent from Theorem~\ref{Théorème de Church-Rosser}.

\end{proof}

\subsubsection{Local Confluence}

We say that $F$ is \emph{locally confluent} if for every $v\in\K{G}$ and for every $T_1,\ T_2\in F$, there exists $v'\in\K{G}$ such that $T_1(v)$ and $T_2(v)$ rewrite into $v'$.

The last result of this section is the formulation of Newman's Lemma \cite{MR0007372} in terms of reduction operators.

\subsubsection{Proposition}\label{Lemme de Newman}

\emph{The set F is confluent if and only if it is locally confluent.}

\begin{proof}

Assume that $F$ is confluent. Let $v$ be an element of $\K{G}$ and let $T_1,\ T_2\in F$. Let $i=1$ or 2. From Theorem~\ref{Théorème de Church-Rosser}, $T_i(v)$ rewrites into $\left(\wedge F\right)\left(T_i(v)\right)$. The latter is equal to $\left(\wedge F\right)(v)$ from Relation~(\ref{Characterisation of the order relation}) (see~\ref{Lattice structure}). Hence, $F$ is locally confluent.

Assume that $F$ is locally confluent. From Proposition~\ref{Confluence vs unique forme normale}, it is sufficient to show that every element $v$ of $\K{G}$ admits a unique $F$-normal form. We show this assertion by induction on $v$. If $v$ is equal to 0, there is nothing to prove. Let $v$ be an element of $\K{G}$. Assume by induction that for every $v'\leq_\text{mul}v$, $v'$ admits a unique $F$-normal form. If $v$ belongs to $\K{\red{F}}$, then $v$ admits a unique $F$-normal form which is itself. Assume that $v$ does not belong to $\K{\red{F}}$. Let $v_1$ and $v_2$ be two $F$-normal forms of $v$. For $i=1$ or $2$, there exists $R_i\in\EV{F}$ such that $v_i$ is equal to $R_i(v)$. We write $R_i=R'_i\circ T_i$, where $T_i$ and $R'_i$ belong to $F$ and $\EV{F}$, respectively. The operator $T_i$ is chosen in such a way that $T_i(v)$ is different from $v$. The set $F$ being locally confluent, there exists $u\in\K{G}$ such that $T_i(v)$ rewrites into $u$. From Proposition~\ref{Existence d'une forme normale}, $u$ admits an $F$-normal form $\widehat{u}$. The latter is also an $F$-normal form of $T_i(v)$. Moreover, $v_i$ is equal to $R'_i\left(T_i(v)\right)$, so that it is also an $F$-normal form of $T_i(v)$. The latter is strictly smaller than $v$ for $\leq_\text{mul}$. By the induction hypothesis, it admits a unique $F$-normal form, so that $v_i$ is equal to $\widehat{u}$ for $i=1$ or $2$. In particular, $v_1$ is equal to $v_2$, so that $v$ admits a unique $F$-normal form.

\end{proof}

\subsection{Reduction Operators and Abstract Rewriting}\label{reduction operators and abstract rewriting}

We fix a subset $F$ of $\RO$.

\subsubsection{Abstract Rewriting Systems}

An \emph{abstract rewriting system} is a pair $\left(A,\F{}\right)$ where $A$ is a set and $\F{}$ is a binary relation on $A$. We write $a\F{}b$ instead of $\left(a,b\right)\in\F{}$. We denote by $\overset{+}{\F{}}$, $\overset{*}{\F{}}$ and $\overset{*}{\underset{}{\longleftrightarrow}}$ the transitive closure, the reflexive transitive closure and the reflexive transitive symmetric closure of $\F{}$, respectively. Finally, we recall the notion of normal form: a normal form of $\left(A,\F{}\right)$ is an element $a$ of $A$ such that there does not exist any element of $\F{}$ with the form $a\F{}b$.

\subsubsection{Confluence and Church-Rosser Property}

Let $\left(A,\F{}\right)$ be an abstract rewriting system. We say that $\F{}$ is \emph{confluent} if for every $a,\ b,\ c\in A$ such that $a\overset{*}{\F{}}b$ and $a\overset{*}{\F{}}c$, there exists $d\in A$ such that $b\overset{*}{\F{}}d$ and $c\overset{*}{\F{}}d$.  We say that $\F{}$ has the \emph{Church-Rosser property} if for every $a,\ b\in A$ such that $a\overset{*}{\underset{}{\longleftrightarrow}} b$, there exists $c\in A$ such that $a\overset{*}{\F{}}c$ and $b\overset{*}{\F{}}c$. Recall from~\cite[Theorem 2.1.5]{MR1629216} that $\F{}$ is confluent if and only if it has the Church-Rosser property.

\subsubsection{Definition}\label{Definition of the abstract rewriting system associated with a set of reduction operators}

We consider the abstract rewriting system $\left(\K{G},\F{F}\right)$ defined by $v\F{F}v'$ if and only if there exists $T\in F$ such that $v$ does not belong to $\K{\red{T}}$ and $v'$ is equal to $T(v)$. 

\subsubsection{Remarks}

\begin{enumerate}
\item If we have $v\F{F}v'$, then we have $v'\leq_\text{mul}v$. The order $\leq_\text{mul}$ being well-founded, the relation $\F{F}$ is also well-founded.
\item We have $v\overset{*}{\F{F}}v'$ if and only if there exists $R\in\EV{F}$ such that $v'$ is equal to $R(v)$, that is, if and only if $v$ rewrites into $v'$. In particular, $v\overset{*}{\F{F}}v'$ implies that $v-v'$ belongs to $\ker\left(\wedge F\right)$.
\end{enumerate}

\subsubsection{Lemma}\label{Linear stability}

\emph{Let $v_1,\ v_2,\ v_3\in\K{G}$ such that $v_1\overset{*}{\underset{F}{\longleftrightarrow}}v_3$. Then, we have $v_1+v_2\overset{*}{\underset{F}{\longleftrightarrow}}v_2+v_3$.}

\begin{proof}

For every $u_1,\ u_2\in\K{G}$ and for every $T\in F$, we have $u_1+u_2\overset{*}{\F{F}}T(u_1+u_2)$ as well as $u_1+T(u_2)\overset{*}{\F{F}}T(u_1+u_2)$. Hence, we have
\begin{equation} \label{Equivalence relation}
u_1+u_2\overset{*}{\underset{F}{\longleftrightarrow}}u_1+T(u_2).
\end{equation}
Let $u_3\in\K{G}$ such that $u_2\overset{*}{\F{F}}u_3$, that is, there exists $R\in\EV{F}$ such that $u_3$ is equal to $R(u_2)$. From Relation~(\ref{Equivalence relation}), we have 
\begin{equation} \label{Equivalence relation 2}
u_1+u_2\overset{*}{\underset{F}{\longleftrightarrow}}u_1+u_3.
\end{equation}
Let $v_1,\ v_2\in\K{G}$ such that $v_1\overset{*}{\underset{F}{\longleftrightarrow}}v_2$, that is, there exists a zig-zag
\[v_1=u_1\overset{*}{\F{F}}u_2\overset{*}{\underset{F}{\longleftarrow}}u_3\overset{*}{\underset{F}{\longrightarrow}}\cdots\overset{*}{\F{F}}u_{r-1}\overset{*}{\underset{F}{\longleftarrow}}u_r=v_2.\]
Relation~(\ref{Equivalence relation 2}) implies that for every $v_3\in\K{G}$ and for every $i\in\{1,\cdots,r-1\}$, we have $u_i+v_3\overset{*}{\underset{F}{\longleftrightarrow}}u_{i+1}+v_3$. Thus, Lemma~\ref{Linear stability} holds.

\end{proof}

\subsubsection{Proposition}\label{Characterisation of the F-reduction system}

\emph{For every $v_1,\ v_2~\in~\K{G}$, we have}
\[v_1\overset{*}{\underset{F}{\longleftrightarrow}}v_2\ if\ and\ only\ if\ v_1-v_2\in\ker\left(\wedge F\right).\]

\begin{proof}

Assume that $v_1\overset{*}{\underset{F}{\longleftrightarrow}}v_2$. Hence, there exists a zig-zag
\[v_1=u_1\overset{*}{\F{F}}u_2\overset{*}{\underset{F}{\longleftarrow}}u_3\overset{*}{\underset{F}{\longrightarrow}}\cdots\overset{*}{\F{F}}u_{r-1}\overset{*}{\underset{F}{\longleftarrow}}u_r=v_2.\]
For every $i\in\{1,\cdots,r-1\}$, $u_i-u_{i+1}$ belongs to $\ker\left(\wedge F\right)$. Hence,
\[v_1-v_2=(u_1-u_2)+(u_2-u_3)+\cdots+(u_{r-2}-u_{r-1})+(u_{r-1}-u_r),\]
belongs to $\ker\left(\wedge F\right)$.

Conversely, assume that $v_1-v_2$ belongs to the kernel of $\wedge F$. The set
\[\left\{v-T(v)\ \mid\ T\in F\ \text{and}\ v\in\K{G}\right\},\]
is a generating set of $\ker\left(\wedge F\right)$. Thus, there exist $T_1, \cdots,T_n\in F$ and $u_1,\cdots,u_n\in\K{G}$ such that $v_1-v_2$ is equal to $\sum_{i=1}^nu_i-T_i(u_i)$, so that we have
\[v_1=\sum_{i=1}^nu_i-T_i(u_i)+v_2.\]
For every $i\in\{1,\cdots,n\}$, we have $u_i-T_i(u_i)\overset{*}{\underset{F}{\longleftrightarrow}}0$. Hence, from Lemma~\ref{Linear stability}, we have $v_1\overset{*}{\underset{F}{\longleftrightarrow}}v_2$.

\end{proof}

\subsubsection{Remark}

Let $v$ be an element of $\K{G}$ and let $[v]$ be the set of elements $v'$ such that $v'-v$ belongs to $\ker\left(\wedge F\right)$. From Proposition~\ref{Characterisation of the F-reduction system}, $[v]$ is the equivalence class of $v$ for the relation $\overset{*}{\underset{F}{\longleftrightarrow}}$. From Lemma~\ref{La borne inf est l'opérateur de normalisation}, $\left(\wedge F\right)(v)$ is the smallest element of this equivalence class.

\subsubsection{Proposition}\label{Church Rosser if and only if Church Rosser}

\emph{The set F has the Church-Rosser property if and only if $\F{F}$ has the Church-Rosser property.}

\begin{proof}

Assume that $F$ has the Church-Rosser property. Let $v,\ v'\in\K{G}$ such that $v\overset{*}{\underset{F}{\longleftrightarrow}}v'$. From Proposition~\ref{Characterisation of the F-reduction system}, $v-v'$ belongs to the kernel of $\wedge F$. We denote by $u$ the common value of $\left(\wedge F\right)(v)$ and $\left(\wedge F\right)(v')$. The set $F$ having the Church-Rosser property, $v$ and $v'$ rewrite into $u$, that is, we have $v\overset{*}{\F{F}}u$ and $v'\overset{*}{\F{F}}u$. Hence, $\F{F}$ has the Church-Rosser property.

Conversely, assume that $\F{F}$ has the Church-Rosser property. Let $v\in\K{G}$. From Proposition~\ref{Characterisation of the F-reduction system}, we have $v\overset{*}{\underset{F}{\longleftrightarrow}}\left(\wedge F\right)(v)$. The relation $\F{F}$ having the Church-Rosser property, there exists $u\in\K{G}$ such that $v\overset{*}{\F{F}}u$ and $\left(\wedge F\right)(v)\overset{*}{\F{F}}u$. Moreover, $\left(wedge F\right)(v)$ belongs to $\K{\red{F}}$, so that it is an $F$-normal form. As a consequence, $u$ is equal to $\left(\wedge F\right)(v)$. Hence, we have $v\overset{*}{\F{F}}\left(\wedge F\right)(v)$, that is, $v$ rewrites into $\left(\wedge F\right)(v)$. That shows that $F$ has the Church-Rosser property.

\end{proof}

\subsubsection{Corollary}\label{Confluence for operators and confluence for abstract rewriting}

\emph{The set F is confluent if and only if $\F{F}$ is confluent.}

\begin{proof}

From Theorem~\ref{Théorème de Church-Rosser}, $F$ is confluent if and only if it has the Church-Rosser property. Hence, from Proposition~\ref{Church Rosser if and only if Church Rosser}, $F$ is confluent if and only if $\F{F}$ has the Church-Rosser property, that is, if and only if $\F{F}$ is confluent.

\end{proof}

\subsubsection{Example}

We consider the pair $P$ of Example~\ref{Illustration of the lattice structure}. We have seen that the pair $P$ is not confluent. The following diagram
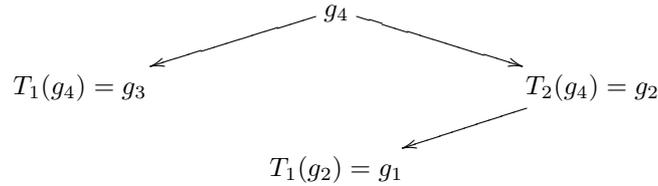
\[
\xymatrix @C = 4em @R = 1.5em{
&
g_4
\ar@1 [rd] 
\ar@1 [ld] 
& \\
T_1(g_4)=g_3
&
&
T_2(g_4)=g_2
\ar@1 [ld]
 \\
&
T_1(g_2)=g_1
& \\
}
\]
shows that we have $g_4\overset{*}{\F{P}}g_1$ and $g_4\overset{*}{\F{P}}g_3$. The two elements $g_1$ and $g_3$ are normal forms, so that $\F{P}$ is not confluent.

\section{Completion and Presentations by Operator}\label{Completion and presentations by operator}

The aim of this section is to formulate algebraically the completion using the lattice structure introduced in Section~\ref{Definition and lattice structure}. We also apply the theory of reduction operators to algebras. Before that, we need to investigate the notion of confluence for a pair of reduction operators.

In Section~\ref{Confluence for a pair of reduction operators title} and Section~\ref{Completion} we fix a well-ordered set $\left(G,<\right)$.

\subsection{Confluence for a Pair of Reduction Operators}\label{Confluence for a pair of reduction operators title}

Throughout this section we fix a pair $P=\left(T_1,T_2\right)$ of reduction operators relative to $\left(G,<\right)$.

\subsubsection{The Braided Products}\label{The braided products}

Given two endomorphisms $S$ and $T$ of $\K{G}$, we denote by $\EV{T,S}^n$ the product $\cdots S\circ T\circ S$ with $n$ factors. Let $g\in G$. From Lemma~\ref{Terminaison}, there exists an integer $n$ such that $\EV{T_2,T_1}^{n}(g)$ and $\EV{T_1,T_2}^{n}(g)$ are $P$-normal forms. Let $n_g$ be the smallest integer satisfying the previous condition. Let $\EV{T_2,T_1}$ and $\EV{T_1,T_2}$ be the two endomorphisms of $\K{G}$ defined by
\[\EV{T_2,T_1}(g)=\EV{T_2,T_1}^{n_g}(g)\ \text{and}\ \EV{T_1,T_2}(g)=\EV{T_1,T_2}^{n_g}(g),\]
for every $g\in G$.

\subsubsection{Remark}\label{Formes normales pour une paire d'opérateurs de réduction}

The vector spaces $\im{\EV{T_2,T_1}}$ and $\im{\EV{T_1,T_2}}$ are included in $\K{\red{P}}$. Hence, every element $v\in\K{G}$ admits at most two $P$-normal forms: $\EV{T_2,T_1}(v)$ and $\EV{T_1,T_2}(v)$.

\subsubsection{Lemma}\label{Confluence for a pair of reduction operators}

\emph{The pair $P$ is confluent if and only if $\EV{T_2,T_1}$ and $\EV{T_1,T_2}$ are equal. In this case we have}
\[\begin{split}
\wedge P&=\EV{T_2,T_1}\\
&=\EV{T_1,T_2}.
\end{split}\]

\begin{proof}
From Proposition~\ref{Confluence vs unique forme normale}, $P$ is confluent if and only if every element of $\K{G}$ admits a unique $P$-normal form. Hence, $P$ is confluent if and only if for every $v\in\K{G}$, $\EV{T_2,T_1}(v)$ and $\EV{T_1,T_2}(v)$ are equal. That shows the first part of the proposition. The second part is a consequence of Lemma~\ref{La borne inf est l'opérateur de normalisation}.
 
\end{proof}

\subsubsection{Dual Braided Products}\label{Dual braided products}

Let $n$ be an integer. We show by induction on $n$ that we have

\begin{equation} \label{Upper-bound}
\begin{split}
\EV{\id{\K{G}}-T_2,\id{\K{G}}-T_1}^n&=\id{\K{G}}+\sum_{i=1}^{n-1}(-1)^i\left(\EV{T_1,T_2}^i+\EV{T_2,T_1}^i\right)+(-1)^{n}\EV{T_2,T_1}^n
,\\
\EV{\id{\K{G}}-T_1,\id{\K{G}}-T_2}^n&=\id{\K{G}}+\sum_{i=1}^{n-1}(-1)^i\left(\EV{T_1,T_2}^i+\EV{T_2,T_1}^i\right)+(-1)^{n}\EV{T_1,T_2}^n.
\end{split}
\end{equation}
We consider the two operators $\EV{\id{\K{G}}-T_1,\id{\K{G}}-T_2}$ and $\EV{\id{\K{G}}-T_2,\id{\K{G}}-T_1}$ defined by
\[\begin{split}
\EV{\id{\K{G}}-T_2,\id{\K{G}}-T_1}(g)&=\EV{\id{\K{G}}-T_2,\id{\K{G}}-T_1}^{n_g}(g),\\
\EV{\id{\K{G}}-T_1,\id{\K{G}}-T_2}(g)&=\EV{\id{\K{G}}-T_1,\id{\K{G}}-T_2}^{n_g}(g),
\end{split}\]
for every $g\in G$.

\subsubsection{Remark}

We deduce from Lemma~\ref{Confluence for a pair of reduction operators} that if the pair $P$ is confluent, then $\EV{\id{\K{G}}-T_2,\id{\K{G}}-T_1}$ and $\EV{\id{\K{G}}-T_1,\id{\K{G}}-T_2}$ are equal. In the sequel, when $P$ is assumed to be confluent, the common value of $\EV{\id{\K{G}}-T_2,\allowbreak{}\id{\K{G}}-T_1}$ and $\EV{\id{\K{G}}-T_1,\allowbreak{}\id{\K{G}}-T_2}$ is denoted by $T$.

\subsubsection{Lemma}\label{Proof of reduction operators}

\emph{Assume that P is confluent. Then, $\idd{\K{G}}-T$ is a reduction operator relative to $\left(G,<\right)$. Moreover, we have}
\[\nred{\id{\K{G}}-T}=\nred{T_1}\cap\nred{T_2}.\]

\begin{proof}

First, we show that $\id{\K{G}}-T$ is a projector. The operators $\id{\K{G}}-T_1$ and $\id{\K{G}}-T_2$ are projectors. Hence, by definition of $T$, for every $g\in G$, and for $i=1$ or 2, we have
\[\left(\id{\K{G}}-T_i\right)\circ T(g)=T(g).\]
Hence, $T$ is a projector, so that $\id{\K{G}}-T$ is also a projector.

Let $g\in G$. Let us show that $g-T(g)$ is either equal to $g$ or strictly smaller than $g$. From Relation~(\ref{Upper-bound}) of~\ref{Dual braided products}, we have
\[g-T(g)=\sum_{i=1}^{n_g-1}(-1)^{i+1}\left(\EV{T_1,T_2}^i+\EV{T_2,T_1}^i\right)(g)+(-1)^{n_g+1}\EV{T_2,T_1}^{n_g}(g).\]
Thus, if $g$ belongs to $\nred{T_1}\cap\nred{T_2}$, then $g-T(g)$ is strictly smaller than $g$. Assume that $g$ does not belong to $\nred{T_1}\cap\nred{T_2}$. Assume that $g$ belongs to $\red{T_1}$ (the case where $g$ belongs to $\red{T_2}$ is analogous). We have
\[\begin{split}
g-T(g)&=g+T_2(g)+\sum_{i=2}^{n_g-1}(-1)^{i+1}\left(\EV{T_1,T_2}^i+\EV{T_1,T_2}^{i-1}\right)(g)+(-1)^{n_g+1}\EV{T_1,T_2}^{n_g-1}(g)\\
&=g.
\end{split}\]
Hence, $\id{\K{G}}-T$ is a reduction operator relative to $\left(G,<\right)$ and the set $\nred{\id{\K{G}}-T}$ is equal to $\nred{T_1}~\cap~\nred{T_2}$.

\end{proof}

\subsubsection{Lemma}\label{Upper-bound for a confluent pair of reduction operators}

\emph{Assume that P is confluent. Then, $T_1\vee T_2$ is equal to $\idd{\K{G}}-T$.}

\begin{proof}

The operator $\id{\K{G}}-T$ being a reduction operator relative to $\left(G,<\right)$ and $\theta$ being a bijection, it is sufficient to show that the kernel of $\id{\K{G}}-T$ equals the one of $T_1\vee T_2$. From Relation~(\ref{Upper-bound}) of~\ref{Dual braided products}, for every $v\in\K{G}$, we have
\[\begin{split}
v-T(v)&=\sum_{i=1}^{n-1}(-1)^{i+1}\left(\EV{T_1,T_2}^i+\EV{T_2,T_1}^i\right)(v)+(-1)^{n+1}\EV{T_2,T_1}^{n}(v)\\
&=\sum_{i=1}^{n-1}(-1)^{i+1}\left(\EV{T_1,T_2}^i+\EV{T_2,T_1}^i\right)(v)+(-1)^{n+1}\EV{T_1,T_2}^{n}(v),
\end{split}\]
where $n$ is an integer greater or equal to $n_g$ for every $g\in G$ belonging to the support of $v$. Hence, $\ker\left(T_1\vee T_2\right)=\ker\left(T_1\right)\cap\ker\left(T_2\right)$ is included in $\ker\left(\id{\K{G}}-T\right)$. Moreover, the operator $\id{\K{G}}-T$ being a projector, its kernel is equal to $\im{T}$, that is, we have
\[\begin{split}
\ker\left(\id{\K{G}}-T\right) &=\im{\EV{\id{\K{G}}-T_2,\id{\K{G}}-T_1}}\\
&=\im{\EV{\id{\K{G}}-T_1,\id{\K{G}}-T_2}}.
\end{split}\]
The latter is included in $\ker\left(T_1\right)$ and $\ker\left(T_2\right)$, so that it is included in $\ker\left(T_1\right)\cap\ker\left(T_2\right)=\ker\left(T_1\vee T_2\right)$.

\end{proof}

\subsubsection{Lemma}\label{Non reduced generators of the upper-bound of a confluent pair}

\emph{Assume that P is confluent. Then, $\nredd{T_1\vee T_2}$ is equal to $\nredd{T_1}~\cap~\nredd{T_2}$.}

\begin{proof}
This is a consequence of Lemma~\ref{Proof of reduction operators} and Lemma~\ref{Upper-bound for a confluent pair of reduction operators}.
\end{proof}

\subsection{Completion}\label{Completion}

We fix a subset $F$ of $\RO$.

\subsubsection{Definitions}\label{Definition of Completion}

\begin{enumerate}
\item A \emph{completion of F} is a subset $F'$ of $\RO$ such that 
\begin{enumerate}
\item $F'$ is confluent,
\item $F\ \subseteq\ F'$ and $\wedge F'=\wedge F$.
\end{enumerate}
\item A \emph{complement of F} is an element $C$ of $\RO$ such that
\begin{enumerate}
\item\label{equivalence relation} $\left(\wedge F\right)\wedge C=\wedge F$,
\item\label{minimality} $\obsred{F}\ \subseteq\ \nred{C}$.
\end{enumerate}
A complement is said to be \emph{minimal} if the inclusion~(\ref{minimality}) is an equality.
\end{enumerate}

\subsubsection{Proposition}\label{Completion and confluence}

\emph{Let} $C\in\RO$ \emph{such that $\left(\wedge F\right)\wedge C$ is equal to $\wedge F$. The set $F\cup\{C\}$ is a completion of F if and only if C is a complement of F.}

\begin{proof}

We let $F'=F\cup\{C\}$. The set $F'$ contains $F$ and is such that $\wedge F'=\wedge F$ by hypothesis. Hence, $F'$ is a completion of $F$ if and only if it is confluent, that is, if and only if $\red{F'}$ is equal to $\red{\wedge F'}$. The set $\red{F'}$ is equal to $\red{F}\cap\red{C}$ and $\wedge F'$ is equal to $\wedge F$. Hence, $F'$ is confluent if and only if we have the following relation
\[\red{F}\cap\red{C}=\red{\wedge F}.\]
Moreover, $\red{F}$ is the disjoint union of $\red{\wedge F}$ and $\obsred{F}$. Hence, we have
\begin{equation}\label{Expression of G F'}
\red{F}\cap\red{C}=\Big(\red{\wedge F}\cap\red{C}\Big)\bigsqcup\left(\obsred{F}\cap\red{C}\right).
\end{equation}
The hypothesis $\left(\wedge F\right)\wedge C$ is equal to $\wedge F$ means that $\wedge F$ is smaller or equal to $C$. Thus, from Lemma~\ref{Inclusion of images}, $\red{\wedge F}$ is included in $\red{C}$. From Relation~(\ref{Expression of G F'}), we have
\[\red{F}\cap\red{C}=\red{\wedge F}\bigsqcup\left(\obsred{F}\cap\red{C}\right).\]
Hence, $F'$ is confluent if and only if $\obsred{F}\cap\red{C}$ is empty, that is, if and only if $C$ is a complement of $F$.

\end{proof}

\subsubsection{Examples}\label{Two examples of completions}

\begin{enumerate}
\item The operator $\wedge F$ is a complement of $F$. In general, this complement is not minimal (see Point~\ref{Second example on Completion}). 
\item\label{Second example on Completion} We consider the pair $P$ of Example~\ref{Illustration of the lattice structure}. Let 
\[C_1=\begin{pmatrix}
1&0&1&0\\
0&1&0&0\\
0&0&0&0\\
0&0&0&1
\end{pmatrix}\ \text{and}\ 
C_2=\begin{pmatrix}
1&0&0&0\\
0&1&1&0\\
0&0&0&0\\
0&0&0&1
\end{pmatrix}.\]
The sets $\nred{C_1}$ and $\nred{C_2}$ are equal to $\{g_3\}$. The latter is equal to $\obsred{P}$ (see Point~\ref{Example of confluence, obstructions} of Example~\ref{Examples for confluence}). Moreover, $\ker\left(C_1\right)$ and $\ker\left(C_2\right)$ are the vector spaces spanned by $g_3-g_1$ and $g_3-g_2$, respectively. These two vector spaces are included in $\ker\left(T_1\wedge T_2\right)$. Hence, $C_1$ and $C_2$ are greater than $T_1\wedge T_2$, that is, we have
\[\left(\wedge P\right)\wedge C_1=\wedge P\ \text{and}\ \left(\wedge P\right)\wedge C_2=\wedge P.\]
We conclude that $C_1$ and $C_2$ are two minimal complements of $P$. We also recall that $\nred{T_1\wedge T_2}$ is equal to $\{g_2,g_3,g_4\}$ (see Example~\ref{Illustration of the lattice structure}), so that $T_1\wedge T_2$ is not a minimal complement of $P$.
\end{enumerate}

\subsubsection{The $F$-Complement}\label{The F-completion}

The \emph{F-complement} is the operator
\[C^F=\left(\wedge F\right)\vee\left(\vee\overline{ F}\right),\]
where $\vee\overline{F}$ is equal to $\theta\left(\K{\red{F}}\right)$.

\subsubsection{Lemma}\label{Confluence of the pair defining the F-completion}

\emph{The pair $P~=~\left(\wedge F,\vee\overline{F}\right)$ is confluent. Moreover, we have}
\begin{equation}\label{equation in lemma}
\nred{\vee\overline{F}}=\red{F}.
\end{equation}

\begin{proof}
Let us show the first part of the lemma. The image of $\wedge F$ is included in $\K{\red{F}}$ which is equal to the kernel of $\vee\overline{F}$. Thus, $\vee\overline{F}\circ\wedge F$ and $\vee\overline{F}\circ\wedge F\circ\vee\overline{F}$ are equal to the zero operator. Hence, we have the equality
\begin{equation}\label{Degree of confluence of the pair defining the F-completion}
\vee\overline{F}\circ\wedge F\circ\vee\overline{F}=\vee\overline{F}\circ\wedge F.
\end{equation}
Hence, from Lemma~\ref{Confluence for a pair of reduction operators}, the pair $P$ is confluent.

Let us show the second part of the lemma. Consider the endomorphism $U$ of $\K{G}$ defined on the basis $G$ in the following way:
\[U(g)=\left\{\begin{split}
& 0,\ \text{if}\ g\in\red{F}\\
& g,\ \text{otherwise}.
\end{split}
\right.\]
The operator $U$ is a projector and is such that for every $g\in G$ such that $U(g)$ is different from $g$, $U(g)$ is equal to $0$. Hence, $U$ is such that for every $g\in G$, we have $U(g)\leq g$, so that it is a reduction operator relative to $\left(G,<\right)$. Moreover, we have
\[\begin{split}
\ker\left(U\right)&=\im{\id{\K{G}}-U}\\
&=\K{\red{F}}.
\end{split}\]
Hence, $U$ and $\vee\overline{F}$ are two reduction operator with same kernel so that they are equal. In particular, $\nred{\vee\overline{F}}$ is equal to $\nred{U}=\red{F}$ which shows the second assertion of the lemma.

\end{proof}

\subsubsection{Theorem}\label{The F-completion is a completion}

\emph{Let F be a subset of} $\RO$. \emph{The F-complement is a minimal complement of F.}

\begin{proof}

By definition, $C^F$ is greater or equal to $\wedge F$, that is, $C^F$ satisfies~(\ref{equivalence relation}) of~\ref{Definition of Completion}.  Let us show that $\obsred{F}$ is equal to $\nred{C^F}$. From Lemma~\ref{Confluence of the pair defining the F-completion}, the pair $\left(\wedge F,\vee\overline{F}\right)$ is confluent. Hence, from Lemma~\ref{Non reduced generators of the upper-bound of a confluent pair} and Relation (\ref{equation in lemma}), we have
\[\begin{split}
\nred{C^F}&=\nred{\left(\wedge F\right)\vee\left(\vee\overline{F}\right)}\\
&=\nred{\wedge F}\cap\nred{\vee\overline{F}}\\
&=\nred{\wedge F}\cap\red{F}\\
&=\obsred{F}.
\end{split}\]

\end{proof}

We end this section with a characterisation of the $F$-complement. For that, we need the following lemma:

\subsubsection{Lemma}\label{Image of the F-completion}

\emph{The set $C^F\left(\obsredd{F}\right)$ is included in} $\K{{\red{\wedge F}}}$.

\begin{proof}

We have seen in Relation~(\ref{Degree of confluence of the pair defining the F-completion}) (see the proof of Lemma~\ref{Confluence of the pair defining the F-completion}) that we have
\[\vee\overline{F}\circ\wedge F\circ\vee\overline{F}=\vee\overline{F}\circ\wedge F.\]
From Lemma~\ref{Confluence of the pair defining the F-completion}, the pair $P~=~\left(\wedge F,\vee\overline{F}\right)$ is confluent. Hence, from Lemma~\ref{Upper-bound for a confluent pair of reduction operators}, we have
\[C^F=\id{\K{G}}-\Big(\id{\K{G}}-\vee\overline{F}\Big)\circ\Big(\id{\K{G}}-\wedge F\Big).\]
The sets $\obsred{F}$ and $\K{\red{\wedge F}}$ are included in $\K{\red{F}}$ which is equal to the kernel of $\vee\overline{F}$. Hence, for every $g\in\obsred{F}$, we have
\[\begin{split}
C^F(g)&=g-\Big(\id{\K{G}}-\vee\overline{F}\Big)\circ\Big(\id{\K{G}}-\wedge F\Big)(g)\\
&=\left(\wedge F\right)(g)+\vee\overline{F}(g)-\left(\vee\overline{F}\circ\wedge F\right)(g)\\
&=\left(\wedge F\right)(g).
\end{split}\]
That shows that $C^F\left(\obsred{F}\right)$ is included in $\K{\red{\wedge F}}$. 

\end{proof}

\subsubsection{Proposition}\label{Characterisation of the F-completion}

\emph{The F-complement is the unique minimal complement C of F such that $C\left(\obsredd{F}\right)$ is included in} $\K{{\red{ F}}}$.

\begin{proof}

Let $C$ be a minimal complement of $F$ such that $C\left(\obsred{F}\right)$ is included in $\K{\red{F}}$. For every $g\in~ G~\setminus~\obsred{F}$, $C^F(g)$ and $C(g)$ are equal to $g$. Thus, it is sufficient to show that for every $g\in\obsred{F}$, $C(g)$ is equal to $C^F(g)$. The set $C\left(\obsred{F}\right)$ being included in $\K{\red{F}}$ and $\nred{C}$ being equal to $\obsred{F}$, $C\left(\obsred{F}\right)$ is in fact included in $\K{\left(\red{F}\setminus\obsred{F}\right)}$, that is, it is included in $\K{\red{\wedge F}}$. From Lemma~\ref{Image of the F-completion}, $C^F\left(\obsred{F}\right)$ is also included in $\K{\red{\wedge F}}$. Hence, for every $g\in\obsred{F}$, we have
\begin{equation}\label{First relation}
\left(\wedge F\circ C^F\right)(g)=C^F(g)\ \text{and}\ \left(\wedge F\circ C\right)(g)=C(g).
\end{equation}
Relation~(\ref{equivalence relation}) of~\ref{Definition of Completion} implies that $C^F$ and $C$ are greater or equal to $\wedge F$. Thus, the equivalence~(\ref{Characterisation of the order relation}) (see~\ref{Lattice structure}) implies that $\wedge F\circ C^F$ and $\wedge F\circ C$ are equal to $\wedge F$. Hence, from Relation~(\ref{First relation}), $C^F(g)$ and $C(g)$ are equal to $\left(\wedge F\right)(g)$ for every $g\in\obsred{F}$, so that Proposition~\ref{Characterisation of the F-completion} holds.

\end{proof}

\subsubsection{Examples}

\begin{enumerate}
\item The operator $C_1$ of Example~\ref{Two examples of completions}, Point~\ref{Second example on Completion} is the $P$-complement.
\item Consider Point~\ref{Non confluence for the braided monoid} of Example~\ref{Examples for confluence}. The $P$-complement maps $yxy$ to $xx$ and fixes all other word.
\end{enumerate}

\subsection{Presentations by Operator}\label{Presentations by operator}

\subsubsection{Algebras}

An associative unitary $\K{}$-algebra is a $\K{}$-vector space \textbf{A} equipped with a $\K{}$-linear map, called \emph{multiplication}, $\mu:\textbf{A}\otimes\textbf{A}\F{}\textbf{A}$ which is associative and for which there exists a unit $1_\textbf{A}$. We say algebra instead of unitary associative $\K{}$-algebra. Given a set $X$, let $X^*$ be the set of words written with $X$. This set admits a monoid structure, where the multiplication is given by concatenation of words and the unit is the empty word. Moreover, the free algebra over $X$ is the vector space $\K{X^*}$ spanned by $X^*$ equipped with the multiplication induced by the one of the monoid $X^*$.

From now on, we fix an algebra \textbf{A}.

\subsubsection{Monomial Orders}

Let $X$ be a set. A \emph{monomial order on $X^*$} is a well-founded total order $<$ on $X^*$ such that the following conditions are fulfilled:
\begin{enumerate}
\item $1<w$ for every word $w$ different from 1,
\item for every $w_1,\ w_2,\ w,\ w'\in X^*$ such that $w<w'$, we have $w_1ww_2<w_1w'w_2$.
\end{enumerate}
In particular, $\left(X^*,<\right)$ is a well-ordered set. In the sequel, given an element $f\in\K{X^*}$, we write $\lm{f}$ (for leading monomial) instead of $\lgen{f}$.

\subsubsection{\G\ Bases}

Let $X$ be a set and let $<$ be a monomial order on $X^*$. Given a subset $E$ of $\K{X^*}$, we let $\lm{E}=\left\{\lm{f}\ \mid\ f\in E\right\}$. Let $I$ be a two-sided ideal of $\K{X^*}$. A subset $R$ of $I$ is called a \emph{\G\ basis of I} if the semi-group ideal spanned by $\lm{R}$ is equal to $\lm{I}$. In other words, $R$ is a \G\ basis of $I$ if and only if for every $w\in\lm{I}$, there exist $w'\in\lm{R}$ and $w_1,\ w_2\in X^*$ such that $w$ is equal to $w_1w'w_2$.

\subsubsection{Definition}

A \emph{presentation by operator} of \textbf{A} is a triple $\EV{\left(X,<\right)\mid S}$, where
\begin{enumerate}
\item $X$ is a set and $<$ is a monomial order on $X^*$,
\item $S$ is a reduction operator relative to $\left(X^*,<\right)$,
\item we have an isomorphism of algebras
\[\textbf{A}\simeq\frac{\K{X^*}}{I\left(\ker\left(S\right)\right)},\]
where $I\left(\ker\left(S\right)\right)$ is the two-sided ideal of $\K{X^*}$ spanned by $\ker\left(S\right)$.
\end{enumerate} 

\subsubsection{Reduction Family of a Presentation}\label{Extensions of S to higher degrees}

Let $X$ be a set and let $n$ be an integer. We denote by $\X{n}$ and $\X{\leq n}$ the set of words of length $n$ and of length smaller or equal to $n$, respectively. Let $\EV{\left(X,<\right)\mid S}$ be a presentation by operator of \textbf{A}. For every integers $n$ and $m$ such that $(n,m)$ is different from $(0,0)$, we let
\[T_{n,m}=\id{\K{\X{\leq n+m-1}}}\ \oplus\ \id{\K{\X{n}}}\otimes S\otimes\id{\K{\X{m}}}.\]
Explicitly, given $w\in X^*$, $T_{n,m}(w)$ is equal to $w$ if the length of $w$ is strictly smaller than $n+m$. If the length of $w$ is greater or equal to $n+m$, we let $w=w_1w_2w_3$, where $w_1$ and $w_3$ have length $n$ and $m$, respectively. Then, $T_{n,m}(w)$ is equal to $w_1S(w_2)w_3$. We also let $T_{0,0}=S$. The \emph{reduction family} of $\EV{\left(X,<\right)\mid S}$ is the set $\left\{T_{n,m},\ 0\leq n,\ m\right\}$.

\subsubsection{Lemma}\label{Extensions of S}

\emph{Let $\EV{\left(X,<\right)\mid S}$ be a presentation by operator of} \textbf{A}. \emph{Let n and m be two integers. Then,} $T_{n,m}$ \emph{is a reduction operator relative to $\left(X^*,<\right)$ and its kernel is equal to} $\K{\X{n}}~\otimes~\ker\left(S\right)~\otimes~\K{\X{m}}$.

\begin{proof}

First, we show that $T_{n,m}$ is a reduction operator relative to $\left(X^*,<\right)$. The operator $S$ being a projector, $T_{n,m}$ is also a projector. Let $w\in X^*$. If the length of $w$ is strictly smaller than $n+m$, then $T_{n,m}(w)$ is equal to $w$. If the length of $w$ is greater or equal $n+m$, we write $w=w_1w_2w_3$, where $w_1$ and $w_3$ have length $n$ and $m$, respectively. If $w_2$ belongs to $\red{S}$, then $T_{n,m}(w)$ is equal to $w$. In the other case, let
\[S(w_2)=\sum\lambda_iw_i,\]
be the decomposition of $S(w_2)$ with respect to the basis $X^*$. We have
\[T_{n,m}(w)=\sum\lambda_iw_1w_iw_3.\]
For every $i\in I$, $w_i$ is strictly smaller than $w_2$. The order $<$ being a monomial order, $w_1w_iw_3$ is strictly smaller than $w$. Hence, $T_{n,m}$ is a reduction operator relative to $\left(X^*,<\right)$. 

Let us show the second part of the lemma. Given an integer $k$, we denote by $\X{\geq k}$ the set of words of length greater or equal to $k$. For every $f\in\K{X^*}$, we write $f=f_1+f_2$, where $f_1$ and $f_2$ are the images of $f$ by the natural projections of $\K{X^*}$ on $\K{\X{\leq n+m-1}}$ and $\K{\X{\geq n+m}}$, respectively. These two vector spaces are stabilised by $T_{n,m}$ and $T_{n,m}(f_1)$ is equal to $f_1$. Thus, $f$ belongs to $\ker\left(T_{n,m}\right)$ if and only if $f_1$ is equal to 0 and $f_2$ belongs to $\ker\left(T_{n,m}\right)$. Moreover, $f_2$ admits a unique decomposition with shape
\[f_2=\sum_{i\in I}w_if_iw'_i,\]
where for every $i\in I$
\begin{enumerate}
\item $w_i$ and $w'_i$ are words of length $n$ and $m$, respectively,
\item for every $j\in I$ such that $j$ is different from $i$, the pair $(w_i,w'_i)$ is different from $(w_j,w'_j)$,
\item $f_i$ is a non zero element of $\K{X^*}$.
\end{enumerate}
We have
\[T_{n,m}\left(f_2\right)=\sum_{i\in I}w_iS(f_i)w'_i.\]
Thus, $f_2$ belongs to $\ker\left(T_{n,m}\right)$ if and only if for every $i\in I$, $f_i$ belongs to $\ker(S)$. Hence, the kernel of $T_{n,m}$ is equal to $\K{\X{n}}\otimes\ker\left(S\right)\otimes\K{\X{m}}$.

\end{proof}

\subsubsection{Remark}

Let $\EV{\left(X,<\right)\mid S}$ be a presentation by operator of \textbf{A}. From Lemma~\ref{Extensions of S}, its reduction family is a subset of $\textbf{RO}\left(X^*,<\right)$. Moreover, the kernel of $\wedge F$ is the sum of the vector spaces $\K{\X{n}}~\otimes~\ker\left(S\right)~\otimes~\K{\X{m}}$, that is, $\ker\left(\wedge F\right)$ is the two-sided ideal spanned by $\ker\left(S\right)$.

\subsubsection{Confluent Presentation}

A \emph{confluent presentation by operator} of \textbf{A} is a presentation by operator of \textbf{A} such that its reduction family is confluent.

\subsubsection{Lemma}\label{Characterisation of the semi-group ideal}

\emph{Let $\EV{\left(X,<\right)\mid S}$ be a presentation by operator of} \textbf{A} \emph{and let F be its reduction family. Let R be the reduced basis of $\ker\left(S\right)$. A word belongs to $\redd{F}$ if and only if it does not belong to the semi-group ideal spanned by $\lmm{R}$.}

\begin{proof}

The set $\red{F}$ is the set of words $w$ such that every sub-word of $w$ belongs to $\red{S}$. From Proposition~\ref{Bijection}, a word belongs to $\red{S}$ if and only if it does not belong to $\lm{R}$, so that Lemma~\ref{Characterisation of the semi-group ideal} holds.

\end{proof}

\subsubsection{Proposition}\label{Characterisation of basis}

\emph{Let $\EV{\left(X,<\right)\mid S}$ be a presentation by operator of} \textbf{A}. \emph{Let R be the reduced basis of $\ker\left(S\right)$. The presentation $\EV{\left(X,<\right)\mid S}$ is confluent if and only if R is a \G\ basis of $I(R)$.}

\begin{proof}

Let $F$ be the reduction family of $\EV{\left(X,<\right)\mid S}$.

Assume that $\EV{\left(X,<\right)\mid S}$ is not confluent. Let $w\in\obsred{F}$. The element $w-\left(\wedge F\right)(w)$ belongs to $I\left(\ker\left(S\right)\right)$. The latter is equal to $I(R)$. Moreover, the leading monomial of $w-\left(\wedge F\right)(w)$ is equal to $w$, so that it belongs to $\red{F}$. Hence, from Lemma~\ref{Characterisation of the semi-group ideal}, $w$ does not belong to the semi-group ideal spanned by $\lm{R}$. Thus, $R$ is not a \G\ basis of $I(R)$.

Assume that $\EV{\left(X,<\right)\mid S}$ is confluent. Let $f\in I(R)$. We assume that $\lc{f}$ is equal to 1. The kernel of $\wedge F$ being equal to $I(R)$, $\left(\wedge F\right)(f)$ is equal to 0. Hence, $\left(\wedge F\right)(\lm{f})$ is equal to $\left(\wedge F\right)(\lm{f}-f)$. The element $\lm{f}-f$ is either equal to $0$ or has a leading monomial strictly smaller than $\lm{f}$. In particular, $\lm{f}$ is not $\wedge F$-reduced. The set $F$ being confluent, that implies that $\lm{f}$ does not belong to $\red{F}$. From Lemma~\ref{Characterisation of the semi-group ideal}, $\lm{f}$ belongs to the semi-group ideal spanned by $\lm{R}$. Thus, $R$ is a \G\ basis of $I(R)$.

\end{proof}

\subsubsection{Theorem}\label{Buchberger algorithm}

\emph{Let $\EV{\left(X,<\right)\mid S}$ be a presentation by operator of \emph{\textbf{A}} and let C be a complement of its reduction family. The triple $\EV{\left(X,<\right)\mid S\wedge C}$ is a confluent presentation of \emph{\textbf{A}}.}

\begin{proof}

We denote by $F$ the reduction family of $\EV{\left(X,<\right)\mid S}$.

First, we show that $\EV{\left(X,<\right)\mid S\wedge C}$ is a presentation of \textbf{A}. For that, we need to show that $I\left(\ker\left( S\wedge C\right)\right)$ is equal to $I\left(\ker\left(S\right)\right)$. The vector space $\ker\left(S\right)$ being included in $\ker\left(S\wedge C\right)$, $I\left(\ker\left(S\right)\right)$ is included in $I\left(\ker\left(S\wedge C\right)\right)$. Moreover, by definition of a complement, $\wedge F$ is smaller or equal to $C$, that is, the kernel of $C$ is included in the one of $\wedge F$. The latter is equal to $I\left(\ker\left(S\right)\right)$. Thus, $\ker\left(S\wedge C\right)$ is also included in $I\left(\ker\left(S\right)\right)$, so that $I\left(\ker\left(S\wedge C\right)\right)$ is included in $I\left(\ker\left(S\right)\right)$.

Let us show that this presentation is confluent. Writing $\tilde{S}=S\wedge C$, we let $\tilde{T}_{0,0}=\tilde{S}$, and for every integers $n$ and $m$ such that $n+m$ is greater or equal to 1, we let 
\[\tilde{T}_{n,m}=\id{\K{\X{\leq n+m-1}}}\ \oplus\ \id{\K{\X{n}}}\otimes\tilde{S}\otimes\id{\K{\X{m}}}.\]
Thus, $\tilde{F}=\left\{\tilde{T}_{n,m},\ 0\leq n,\ m\right\}$ is the reduction family of $\EV{\left(X,<\right)\mid S\wedge C}$. The kernel of $\wedge\tilde{F}$ is equal to $I\left(\ker\left(S\wedge C\right)\right)$. We have seen that the latter is equal to $I\left(\ker\left(S\right)\right)$ which is equal to $\ker\left(\wedge F\right)$. Hence, $\wedge\tilde{F}$ is equal to $\wedge F$. As a consequence, we have to show that $\red{\tilde{F}}$ is included in $\red{\wedge F}$. For every integers $n$ and $m$, $T_{n,m}$ is greater or equal to $\tilde{T}_{n,m}$. Hence, from Lemma~\ref{Inclusion of images}, $\red{\tilde{T}_{n,m}}$ is included in $\red{T_{n,m}}$, so that $\red{\tilde{F}}$ is included in $\red{F}$. Moreover, $\red{\tilde{F}}$ is also included in $\red{S\wedge C}$ which is itself included in $\red{C}$. Hence, $\red{\tilde{F}}$ is included in $\red{F}\cap\red{C}$. From Proposition~\ref{Completion and confluence}, the set $F\cup\{C\}$ is confluent, so that $\red{F}\cap\red{C}$ is equal to $\red{\wedge F}$. That shows that the presentation $\EV{\left(X,<\right)\mid S\wedge C}$ is confluent.

\end{proof}

\subsubsection{Remark}\label{Remark on Buchberger algorithm}

Theorem~\ref{Buchberger algorithm} provides a theoretical method to construct \G\ bases. However, this method is not algorithmic, a priori. Indeed, the operator $C$ which appears in this theorem is relative to the infinite set $X^*$ whereas our implementation~\footnote{http://pastebin.com/0YZCfAD4} of reduction operators requires to work with finite sets. We discuss this problem in the conclusion of the paper.

\section{Generalised Reduction Operators}\label{Generalised reduction operators}

So far, we have been studying reduction operators relative to a well-ordered set. In this section, we investigate the more general case where we do not consider a total order. The purpose of this section is not to provide new results but to explain why the requirement of a total order is of crucial importance in Section~\ref{Rewriting properties of reduction operators} and Section~\ref{Completion and presentations by operator}.

We fix an ordered set $\left(G,<\right)$. 

\subsection{Algebraic Structure}\label{Algebraic structure}

The general definition of reduction operator is stated as follows:

\subsubsection{Definition}\label{General definition of reduction operator}

A \emph{reduction operator relative to} $\left(G,<\right)$ is an idempotent endomorphism $T$ of $\K{G}$ such that for every $g\in G$, we have either $T(g)= g$, or for every $g'$ occurring in the support of $T(g)$, we have $g'<g$. As in the case of well-ordered sets, the set of reduction operators relative to $\left(G,<\right)$ is denoted by $\textbf{RO}\left(G,<\right)$ . Given a reduction operator $T$, the set of $T$-reduced generators is also denoted by $\red{T}$ and its complement in $G$ is denoted by $\nred{T}$.

\subsubsection{From Projectors to Reduction Operators}

As an application, we want to consider a set $G$ together with a non-empty subset $F$ of the set of all linear idempotent endomorphisms of $\K{G}$. We want to equip $G$ with an order $<$ making $F$ a subset of $\RO$. For that, consider the binary relation $<_F$ on $G$ defined by $g'<_Fg$ if there exists $T\in F$ such that $T(g)$ is different from $g$ and such that $g'$ belongs to the support of $T(g)$. The transitive closure of $<_F$ is still denoted by $<_F$. This relation is not necessarily anti-symmetric. Indeed, let $G=\{g_1,g_2\}$ and consider $F=\left(T_1,T_2\right)$, where $T_1$ and $T_2$ are defined by $T_1(g_2)=g_1$, $T_1(g_1)=g_1$, $T_2(g_1)=g_2$ and $T_2(g_2)=g_2$, respectively. Then, we have $g_1<_Fg_2$ and $g_2<_Fg_1$. However, if $<_F$ is well-founded, then it is an order relation, and in this case, $F$ is a subset of $\textbf{RO}\left(G,<_F\right)$. 

\subsubsection{Absence of Reduced Basis}\label{Reduced basis for generalised reduction operators}

We would like to equip the set $\RO$ with a lattice structure. We cannot use the argument of Section~\ref{Definition and lattice structure} because a subspace of $\K{G}$ does not necessarily admit a reduced basis. Indeed, consider $G=\{g_1,g_2,g_3\}$ ordered such: $g_1<g_3$ and $g_2<g_3$. The subspace of $\K{G}$ spanned by $g_3-g_1$ and $g_3-g_2$ does not admit any reduced basis.

In order to equip $\RO$ with an order relation, we need the following lemma:

\subsubsection{Lemma}\label{Inclusion of images, general case}

\emph{Let $T_1$ and $T_2$ be two reduction operators relative to $\left(G,<\right)$ such that $\ker\left(T_1\right)$ is included in $\ker\left(T_2\right)$. Then, $\redd{T_2}$ is included in $\redd{T_1}$.}

\begin{proof}

Assume by way of contradiction that there exists $g\in\red{T_2}$ not belonging to $\red{T_1}$. The element $g-T_1(g)$ belongs to the kernel of $T_1$, so that it belongs to the one of $T_2$. Hence, $T_2(g)$ is equal to $T_2\left(T_1(g)\right)$. The generator $g$ belongs to $\red{T_2}$, so that $T_2(g)$ is equal to $g$. Moreover, $g$ being not $T_1$-reduced, every generator appearing in the support of $T_1(g)$ is strictly smaller than $g$, so that every generator belonging to the support of $T_2\left(T_1(g)\right)$ is also strictly smaller than $g$. Thus, we reach a contradiction.

\end{proof}

\subsubsection{Order Relation}\label{Order relation for generalised reduction operators}

The binary relation defined by $T_1\preceq T_2$ if $\ker\left(T_2\right)\ \subseteq\ \ker\left(T_1\right)$ is clearly reflexive and transitive. Moreover, from Lemma~\ref{Inclusion of images, general case}, if two reduction operators have the same kernel, then they have the same image, so that they are equal. Hence, $\preceq$ is anti-symmetric, so that it is an order relation on $\RO$.

\subsubsection{Absence of a Lattice Structure}\label{Absence of lattice structure}

The order introduced in~\ref{Order relation for generalised reduction operators} does not induce a lattice structure. Consider $G=\{g_1,g_2,g_3,g_4,g_5\}$ ordered such: $g_1<g_3$, $g_1<g_4$, $g_2<g_3$, $g_2<g_4$, $g_3<g_5$ and $g_4<g_5$. Let $T_1,$ $T_2$ be the two reduction operators defined by $\red{T_i}=\{g_1,g_2,g_3,g_4\}$ for $i=1$ or 2, $T_1(g_5)=g_3$ and $T_2(g_5)=g_4$. Consider the two reductions operators $U_1$ and $U_2$ defined by $\red{U_i}~=~\{g_1,g_2\}$ for $i=1$ or 2, $U_1(g_j)~=~g_2$ and $U_2(g_j)=g_1$ for $j\in\{3,4,5\}$. The vector space $\ker\left(T_1\right)+\ker\left(T_2\right)=\K{\{g_5-g_4\}}~\oplus~\K{\{g_5-g_3\}}$ is included in $\ker\left(U_i\right)$ for $i=1$ or 2, that is, $U_i$ is smaller than $T_1$ and $T_2$. Moreover, there does not exist a reduction operator with kernel $\K{\{g_5-g_4\}}~\oplus~\K{\{g_5-g_3\}}$, so that $U_1$ and $U_2$ are two maximal elements smaller than $T_1$ and $T_2$. Hence, $T_1$ and $T_2$ admit a lower bound but they do not admit a greatest lower bound. Moreover, even when a greatest lower bound exists, its kernel is not necessarily the sum of the kernels. Consider the example from~\ref{Reduced basis for generalised reduction operators}: $G=\{g_1,g_2,g_3\}$ with $g_1<g_3$ and $g_2<g_3$, and let $P=\left(T_1,T_2\right)$ where, for $i=1$ or 2, $\nred{T_i}=\{g_3\}$ and $T_i(g_3)=g_i$. Then, we check that $T_1$ and $T_2$ admit a lower bound which is the zero operator, that is, the kernel of this lower bound is equal to $\K{G}$.

\subsection{Rewriting Properties}\label{Rewriting properties, general case}

In this section, we investigate the rewriting properties associated to generalised reduction operators. We have seen in the previous section that, given a subset $F$ or $\RO$, there does not necessarily exist a reduction operator with kernel $\sum_{T\in F}\ker\left(T\right)$. Hence, in order to define the notion of confluence as it was done in~\ref{Confluence of a finite family of reduction operators}, we have to consider subsets of $\RO$ for which such a reduction operator exists. For that, we introduce the following definition:

\subsubsection{Completable Sets}

A subset $F$ of $\RO$ for which there exists a reduction operator with kernel $\sum_{T\in F}\ker\left(T\right)$ is said to be \emph{completable}. 

\subsubsection{Confluence and Church-Rosser Property}

Let $F$ be a completable set. The reduction operator whose kernel is equal to $\sum_{T\in F}\ker\left(T\right)$ is denoted by $\wedge F$. From Lemma~\ref{Inclusion of images, general case}, $\red{\wedge F}$ is included in $\red{T}$, so that the set $\obsred{F}$ is well-defined. The set $F$ is said to be \emph{confluent} if $\obsred{F}$ is the empty set. We say that $v$ \emph{rewrites into} $v'$ as it was done in Section~\ref{F-compositions} and that $F$ has the \emph{Church-Rosser property} if for every $v\in\K{G}$, $v$ rewrites into $\left(\wedge F\right)(v)$. Finally, the binary relation $\F{F}$ on $\K{G}$ is defined as it was done in Section~\ref{reduction operators and abstract rewriting}.

\subsubsection{Normalising Relations}

In Theorem~\ref{Generalised confluence}, we use the notion of \emph{normalising relation.} Let us recall it. Given an abstract rewriting system $\left(A,\F{}\right)$, the relation $\F{}$ is said to be normalising if every element of $A$ admits at least one normal form.

\subsubsection{Theorem}\label{Generalised confluence}

\emph{Let F be a completable subset of} $\RO$. \emph{The following assertions are equivalent:}
\begin{enumerate}
\item\label{generalised confluence} \emph{F is confluent and $\F{F}$ is normalising,}
\item\label{generalised C-R} \emph{F has the Church-Rosser property,}
\item\label{confluence of the relation for generalised RO} \emph{$\F{F}$ is confluent.}
\end{enumerate}

\begin{proof}

The proof of the equivalence between~\ref{generalised C-R} and~\ref{confluence of the relation for generalised RO} is the same as in Proposition~\ref{Church Rosser if and only if Church Rosser} (indeed, we check that in Section~\ref{reduction operators and abstract rewriting} we only require the existence of the operator $\wedge F$). As in Theorem~\ref{Théorème de Church-Rosser}, we show that~\ref{generalised C-R} implies that $F$ is confluent. Moreover, if~\ref{generalised C-R} holds, every element $v$ of $\K{G}$ rewrites into $\left(\wedge F\right)(v)$, that is, we have $v\overset{*}{\F{F}}\left(\wedge F\right)(v)$. The latter belongs to $\K{\red{F}}$, so that it is a normal form for $\F{F}$. Hence, $\F{F}$ is normalising. Thus,~\ref{generalised C-R} implies~\ref{generalised confluence}. Assume~\ref{generalised confluence} and let us show~\ref{confluence of the relation for generalised RO}. Let $v_1,\ v_2,\ v_3\in\K{G}$ such that $v_1\overset{*}{\F{F}}v_2$ and $v_1\overset{*}{\F{F}}v_3$. The relation $\F{F}$ being normalising, $v_2$ and $v_3$ admit normal forms. Let $\widehat{v_2}$ and $\widehat{v_3}$ be normal forms of $v_2$ and $v_3$, respectively. We have $\widehat{v_2}\overset{*}{\underset{F}{\longleftrightarrow}}\widehat{v_3}$, so that $\widehat{v_2}-\widehat{v_3}$ belongs to the kernel of $\wedge F$. Moreover, $\widehat{v_2}$ and $\widehat{v_3}$ being normal forms, they belong to $\K{\red{F}}$, that is, $\K{\red{\wedge F}}$ since $F$ is confluent. Hence, $\widehat{v_2}-\widehat{v_3}$ also belongs to $\K{\red{\wedge F}}$, that is, it belongs to the image of $\wedge F$. Thus, $\widehat{v_2}-\widehat{v_3}$ belongs to $\ker\left(\wedge F\right)\cap\im{\wedge F}$ which is reduced to $\{0\}$ since $\wedge F$ is a projector. We conclude that $\widehat{v_2}$ is equal to $\widehat{v_3}$, so that $\F{F}$ is confluent.

\end{proof}

\subsubsection{Completion}

Given a completable set $F$, the notion of \emph{complement} is stated as follows: a complement of $F$ is a reduction operator $C$ satisfying
\begin{enumerate}
\item\label{completion is greater} $\wedge F\preceq C$,
\item $\obsred{F}\ \subseteq\ \nred{C}$.
\end{enumerate}

\subsubsection{Completion and Rewriting}

Let $C$ be a reduction operator, greater or equal to $\wedge F$. We write $F'=F\cup\{C\}$. The vector space $\sum_{T\in F'}\ker\left(T\right)$ is equal to $\sum_{T\in F}\ker\left(T\right)$. Hence, the set $F'$ is also completable and $\wedge F'$ is equal to $\wedge F$. In particular, the two equivalence relations $\overset{*}{\underset{F}{\longleftrightarrow}}$ and $\overset{*}{\underset{F'}{\longleftrightarrow}}$ are equal. Using the arguments of the proof of Proposition~\ref{Completion and confluence}, we get

\subsubsection{Proposition}

\emph{Let} $C\in\RO$ \emph{such that C is greater or equal to $\wedge F$. The set $F\cup\{C\}$ is confluent if and only if C is a complement of F.}

\subsubsection{Remark}

We have seen that the absence of a total order on $G$ implies that the reduction operator $\wedge F$ does not necessarily exist. In particular, the notion of complement is not necessarily defined, like the Knuth-Bendix completion algorithm in term rewriting does not necessarily succeed. Consider the example from~\ref{Reduced basis for generalised reduction operators}: $G=\{g_1,g_2,g_3\}$ with $g_1<g_3$ and $g_2<g_3$, and let $P~=~\left(T_1,T_2\right)$ where, for $i=1$ or 2, $\nred{T_i}=\{g_3\}$ and $T_i(g_3)=g_i$. We have seen in~\ref{Reduced basis for generalised reduction operators} that $T_1\wedge T_2$ does not exist. Moreover, the Knuth-Bendix completion algorithm does not work because as shown on the following diagram
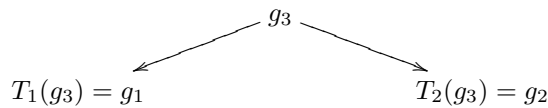
\[
\xymatrix @C = 4em @R = 1.5em{
&
g_3
\ar@1 [rd] 
\ar@1 [ld] 
& \\
T_1(g_3)=g_1
&
&
T_2(g_3)=g_2
}
\]
$g_3$ admits two distinct normal forms $g_1$ and $g_2$, and these two normal forms cannot be compared.\newline

\paragraph{Conclusion.}

The approach using reduction operators shows that some familiar concepts of rewriting, such as confluence or completion, can be formulated in terms of algebraic conditions. In particular, we related the algebraic formulation of confluence to \G\ bases. Even if reduction operators provides a theoretical method to construct \G\ bases as stated in Theorem~\ref{Buchberger algorithm}, we have seen in Remark~\ref{Remark on Buchberger algorithm} that this method is not algorithmic. A natural extension of our work is to transform this method into a procedure. For that, a local criterion has to be used to restrict ourselves to finite-dimensional subspaces of $\K{X^*}$. Such criterion exists for \G\ bases, using \emph{S-elements}~\cite[Section 5.3]{MR1299371}. A future work consists in adapting this criterion to reduction operators in order to obtain an effective procedure.

\bibliography{Biblio}

\newpage\noindent
\textbf{Cyrille Chenavier}\\
\textbf{INRIA, \'equipe $\pi r^2$}\\
\textbf{Laboratoire IRIF, CNRS UMR 8243}\\
\textbf{Universit\'e Paris-Diderot}\\
\textbf{Case 7014}\\
\textbf{75205 PARIS Cedex 13}\\
\textbf{chenavier@pps.univ-paris-diderot.fr}

\end{document}